\documentclass[a4paper,11pt]{article}

\usepackage{amsfonts,amsmath,amssymb,bm}
\usepackage{color,cases,float}
\usepackage{graphicx}
\usepackage{hyperref}
\usepackage{psfrag}
\usepackage{eucal}
\usepackage{algorithm}
\usepackage{overpic}
\usepackage{subfigure}
\usepackage{url}
\usepackage[all,cmtip]{xy}
\usepackage{algorithm}
\usepackage[noend]{algpseudocode}

\definecolor{darkblue}{rgb}{0.0,0.0,0.6}
\hypersetup{colorlinks,breaklinks,
	linkcolor=darkblue,urlcolor=darkblue,
	anchorcolor=darkblue,citecolor=darkblue}
\newtheorem{assumption}{Assumption}
\newtheorem{definition}{Definition}
\newtheorem{lem}{Lemma}
\newtheorem{rem}{Remark}
\newenvironment{proof}{\paragraph{Proof:}}{\hfill$\square$}
\newtheorem{prop}{Proposition}
\newtheorem{theorem}{Theorem}

\newcommand{\E}{\mathrm{E}}

\newcommand{\R}{\mathrm{Re}}
\newcommand{\Proj}{\mathrm{Proj}}

\def\R{\mathbb{R}}

\def\zb{\mathbf{z}}
\def\bx{\boldsymbol{x}}

\def\ab{\boldsymbol{a}}

\def\mb{\boldsymbol{m}}

\def\yb{\mathbf{y}}
\def\Prt{\mbox{Proj}_{(1-\rho_t)\boldsymbol{A}}}
\def\s{\sigma}
\def\g{\gamma}
\def\by{\mathbf{y}}
\def\bmu{\boldsymbol{\mu}}

\newcommand{\tat}[1]{\normalsize{{\color{red}\ #1}}}

\def\tat#1{{\color{red}#1}}

\def\R{\mathbb{R}}

\def\zb{\mathbf{z}}
\def\xb{\boldsymbol{x}}

\def\ba{\boldsymbol{a}}

\def\bmu{\boldsymbol{\mu}}

\def\Mb{\boldsymbol{M}}
\def\Qb{\boldsymbol{Q}}

\def\Gb{\boldsymbol{G}}

\def\Hb{\boldsymbol{H}}
\def\Rb{\boldsymbol{R}}
\def\Ab{\boldsymbol{A}}
\def\Pb{\boldsymbol{P}}

\def\bxi{\boldsymbol{\xi}}
\def\r{~}

\usepackage{anysize}
\marginsize{2cm}{2cm}{1.5cm}{1.5cm}

\title{On the Rate of Convergence of Payoff-based Algorithms to Nash Equilibrium in  Strongly Monotone Games}
\author{Tatiana Tatarenko, \href{mailto:tatiana.tatarenko@rmr.tu-darmstadt.de}{tatiana.tatarenko@rmr.tu-darmstadt.de} \\
       Maryam Kamgarpour, \href{mailto:maryam.kamgarpour@epfl.ch}{maryam.kamgarpour@epfl.ch} }

\begin{document}
\maketitle

\begin{abstract}
We derive the rate of convergence to Nash equilibria for  the payoff-based algorithm proposed in \cite{tat_kam_TAC}. These rates are achieved under the standard assumption of convexity of the game, strong monotonicity and differentiability of the pseudo-gradient. In particular, we show the algorithm achieves  $O(\frac{1}{T})$ in the two-point function evaluating setting and $O(\frac{1}{\sqrt{T}})$ in the one-point function evaluation under additional requirement of Lipschitz continuity of the pseudo-gradient. These rates are to our knowledge the best known rates for the corresponding problem classes.
\end{abstract}

\section{Introduction}
\label{sec:intro}
% what's the problem about
The problem of payoff-based learning  of Nash equilibrium in a multi-agent game concerns deriving a learning algorithm for each player  that uses only evaluation of the player's cost function. The payoff-based information in a game corresponds to the so-called zero-order oracle setting in optimization. Payoff-based algorithms for games over continuous action spaces have been independently proposed in \cite{bravo2018bandit,tat_kam_TAC}. The underlying assumptions have been the convexity of the game and the strong monotonicity of the game pseudo-gradient. While both works derived the rate of convergence to a Nash equilibrium, the lower bound for the convergence rate in this problem class were unknown. Thus, an open question was whether one can achieve faster convergence rates for the above problem class.

% existing rates for payoff-based
Letting $T$ denote the number of iterations of the algorithm, it was shown in \cite{bravo2018bandit} that the proposed algorithm will ensure convergence to the Nash equilibrium at a rate of $O(\frac{1}{T^{1/3}})$. Though an explicit rate was not derived in \cite{tat_kam_TAC}, it could be verified from Theorem 3 in that work that the proposed algorithm under the optimal  parameters would also achieve a rate of $O(\frac{1}{T^{1/3}})$. In this work, we show that the algorithm proposed in \cite{tat_kam_TAC} achieves a rate of $O(\frac{1}{T^{1/2}})$. This improved rate is derived with a more refined analysis of the algorithm rather than modifying it or its parameters.  Furthermore, we propose a modified algorithm to account for  the possibility of two-point evaluation of the players' cost function. For this case, we show that the convergence rate can be improved to $O(\frac{1}{T})$. We expect these rates to be optimal as we will reason below.

% optimal rates for payoff-based
Convex games with monotone pseudo-gradient include the class of convex optimization problems. In particular,  if the game is potential, which is the case if the pseudo-gradient has symmetric Jacobian, then computing Nash equilibrium is equivalent to computing a minimum of a convex potential function. Now, strongly convex smooth optimization corresponds to strongly monotone games with Lipschitz pseudo-gradients. Hence, the lower bound for this class of optimizations would imply a lower bound for the corresponding class of games. It is known that the lower bound on the zero-order  strongly convex smooth optimization is $O(\frac{1}{T^{1/2}})$ in case of one-point function evaluation \cite{shamir2013complexity} and $O(\frac{1}{T})$ in case of full gradient information \cite{Rakhlin2012MakingGD}. Since we show that our algorithm achieves these rates under the payoff-based information, it appears that we achieve the optimal rate in terms of iteration complexity. However, in the one-point setting, we have an additional assumption on the differentiability of the cost functions' gradients (see Assumptions~\ref{assum:Lipschitz}) in comparison to smoothness in \cite{shamir2013complexity} (see Section 2 therein). This assumption is required to analyze strongly monotone properties of the pseudo-gradient in the mixed strategies (see Proposition~\ref{prop:strMon}). Whether a better lower bound can be achieved for the class of optimization problems and games under  this additional assumption is an open question.

% optimal rate for first-order info
As for our algorithmic approach and the derived rates few remarks are in order. First, to estimate the gradient of a player's cost function with respect to her decision variable,  using payoff-based  information, we apply a randomized sampling strategy that estimates the gradient of a smoothed version of the cost function. The randomization approach is similar in spirit to \cite{flaxman2005online} for bandit online learning but in contrast to the above work, we use a Gaussian distribution  inspired by stochastic sampling in \cite{Thatha, nesterov2011random}. The introduced noise term has a diminishing variance denoted by $\sigma_t$. We can then show that in the two-point evaluation, the variance of the estimated gradient is bounded. On the other hand, under the one-point evaluation, the variance is of order $O(\frac{1}{\sigma_t^2})$. This limits the rate at which $\sigma_t$ can converge to zero to ensure convergence to the Nash equilibrium in the one-point setting.

% comparison
We also note that the proposed payoff-based algorithm in the one-point evaluation setting, apart from the sampling distribution, can be considered as a subclass of the  mirror-descent approach in \cite{bravo2018bandit} (by using the Euclidean distance as the  regularizer). The choice of Gaussian distribution requires a different approach to ensure feasibility of the algorithm iterates and different proof techniques. However, this choice is not the reason for the improved rates achieved here. These improvements are due to the fact that unlike the analysis in \cite{bravo2018bandit, tat_kam_TAC}, in this work we compare the algorithm iterates to the Nash equilibria of a corresponding time-varying smoothed game. This modified analysis technique is inspired by \cite{PerchetBach}, who considered zero-order optimization under different smoothness assumptions. 

% contributions
Our contributions are thus as follows. We propose an algorithm to learn Nash equilibria in convex games with strongly monotone differentiable pseudo-gradient under one-point and two-point function evaluations. We establish its rate of convergence as $O(\frac{1}{\sqrt{T}})$ (Theorem \ref{th:main}) in the two-point setting and as $O(\frac{1}{T})$ in the one-point zero-order setting (Theorem \ref{th:main2}). To derive these rates, in the one-point setting, we introduce a new analysis technique in which the iterates of the algorithm are compared with the Nash equilibria of a time-varying smooth game  whose pseudo-gradient is defined in \eqref{eq:smooth_t_game}. The supporting results are in  Propositions \eqref{prop:tildeM}--\eqref{prop:y-mu}.  In the two-point setting, we compare the algorithm iterates to the iterates arising from a modified procedure  defined in \eqref{eq:z} based on exact pseudo-gradient information. Our  supporting results in Propositions \eqref{prop:diffz}-\eqref{prop:z-mu} lead to the rate of convergence reported.

% organization
The rest of the paper is organized as follows. In Section \ref{sec:problem} we formulate the problem of payoff-based learning and state the assumptions on the considered class of games. In Section \ref{sec:procedure} we detail the proposed payoff-based approaches in the one-point and two-point setting. In Section~\ref{sec:rate} we state the main theorems and proofs on the convergence rates of the algorithms. The proofs of all supporting lemmas and propositions are provided in the Appendix. In Section \ref{sec:conclusion} we conclude the paper.

\section{Problem Formulation}
\label{sec:problem}
\allowdisplaybreaks
Consider a game $\Gamma (N, \{A_i\}, \{J_i\})$ with $N$ players, the sets of players' actions $A_i\subseteq \R^d$, $i\in[N]$, and the cost (objective) functions $J_i:\Ab\to\R$, where $\Ab = A_1\times\ldots\times A_N$ denotes the set of joint actions. We restrict the class of games as follows.
%We assume the following regarding  $\Gamma (N, \{A_i\}, \{J_i\})$.

\begin{assumption}\label{assum:convex}
	The game under consideration is \emph{convex}. Namely, for all $i\in[N]$ the set $A_i$ is convex and closed, the cost function $J_i(\ba^i, \ba^{-i})$ is defined on $\R^{Nd}$, continuously differentiable in $\ba$ and convex in $\ba^i$ for  fixed $\ba^{-i}$.
\end{assumption}

\begin{assumption}\label{assum:compact}
	The action sets $A_i$ are compact for all $i\in[N]$.
	%than a polynomial function as $\|\ba\|\to\infty$.
\end{assumption}

\begin{assumption}\label{assum:CG_grad}
	The mapping $\Mb:\R^{Nd}\to\R^{Nd}$, referred to as the \emph{pseudo-gradient}, defined by
	\begin{align}\label{eq:gamemapping}
		\nonumber
		\Mb(\ba) &= (\nabla_{\ba^i} J_i(\ba^i, \ba^{-i}))_{i=1}^N=(\Mb_1(\ba), \ldots, \Mb_N(\ba))^{\top},\\ \nonumber
		&\mbox{where }\Mb_i(\ba) = (M_{i,1}(\ba), \ldots, M_{i,d}(\ba))^{\top}, \;
		M_{i,k}(\ba)&= \frac{\partial J_i(\ba)}{\partial a^i_k},\quad \ba\in\Ab, \quad i\in[N], \quad k\in[d],
	\end{align}
	is \emph{strongly monotone on $\R^{Nd}$} with the constant $\nu$.
\end{assumption}

We consider a \emph{Nash equilibrium} in game $\Gamma (N, \{A_i\}, \{J_i\})$  as a stable solution outcome because it represents a joint action from which no player has any incentive to unilaterally deviate.

\begin{definition}\label{def:NE}
	A point $\ba^*\in\Ab$ is called a \emph{Nash equilibrium} if for any $i\in[N]$ and $\ba^i\in A_i$
	\[J_i(\ba^{i*},\ba^{-i*})\le J_i(\ba^{i},\ba^{-i*}).\]
\end{definition}

Our goal is to learn such a stable action in a game through designing a payoff-based algorithm. To do so, we first connect existence of Nash equilibria for $\Gamma (N, \{A_i\}, \{J_i\})$ with solution set of a corresponding variational inequality problem.
\begin{definition}
	Consider a mapping $\boldsymbol T(\cdot)$: $\R^d \to \R^d$ and a set $Y \subseteq \R^d$. The
	\emph{solution set $SOL(Y,\boldsymbol T)$ to the variational inequality problem} $VI(Y,\boldsymbol T)$ is the set of vectors $\by^* \in Y$ such that $(\boldsymbol T(\by^*), \by-\by^*) \ge 0$, $\forall \by \in Y$.
\end{definition}

\begin{theorem}\label{th:VINE}(\cite[Proposition\r1.4.2]{FaccPang1})
	Given a game $\Gamma (N, \{A_i\}, \{J_i\})$ with pseudo-gradient $\Mb$, suppose that the action sets $\{A_i\}$ are closed and convex, the cost functions $\{J_i\}$ are continuously differentiable in $\ba$ and convex in $\ba^i$ for every fixed $\ba^{-i}$ on the interior of $\Ab$. Then, some vector $\ba^*\in \Ab$ is a Nash equilibrium in $\Gamma$, if and only if $\ba^*\in SOL(\Ab,\boldsymbol M)$. Moreover, if additionally Assumptions~\ref{assum:compact} holds, the Nash equilibrium in $\Gamma$ exists. If Assumptions~\ref{assum:CG_grad} holds, the Nash equilibrium in $\Gamma$ exists and is unique.
\end{theorem}

%\begin{rem}\label{rem:inftybeh}
%  If additionally to Assumption~\ref{assum:Lipschitz} for each $i\in[N]$ the condition
%  \[\|\nabla^2J_i({\bx_1}) - \nabla^2 J_i(\bx_2)\|\le K_1\|\bx_1-\bx_2\| \]
%  holds for all $\bx_1$,$\bx_2\in\R^{Nd}$ with some $K_1$, then for all $i \in [N]$, $J_i(\bx) = O(\|\bx\|^3)$ as  $\|\bx\|\to\infty$. Indeed, according to the Taylor's expansion, for any $\bx\in\R^{Nd}$ and $\bmu\in\Ab$
%  \begin{align*}
	%    J_i(\bx)  = &J_i(\bmu) + (\nabla J_i(\bmu), \bx-\bmu) + ((\bx-\bmu)^T\nabla^2 J_i(\tilde{\bx}),\bx-\bmu)\\
	%     =&J_i(\bmu) + (\nabla J_i(\bmu), \bx-\bmu) + ((\bx-\bmu)^T\nabla^2 J_i(\bmu),\bx-\bmu)\\
	%     & + ((\bx-\bmu)^T[\nabla^2 J_i(\tilde{\bx}) - \nabla^2 J_i(\bmu)],\bx-\bmu),
	%  \end{align*}
%  where $\tilde{\bx} = \bx + \theta(\bmu-\bx)$ for some $\theta\in[0,1]$. Since $\bmu$
% is from the compact set and by taking into account Assumption~\ref{assum:Lipschitz}, we obtain
%   \begin{align*}
	%    |J_i(\bx)| &\le  O(1) + O(1) \|\bx-\bmu\|  + O(1)\|\bx-\bmu\|^2 +\|\bx-\bmu\|^2 K_1\|\tilde{\bx}-\bmu\| \cr
	%      & = O(\|\bx\|^3), \, \mbox{as $\|\bx\|\to\infty$}.
	%  \end{align*}
%  Thus, in  this case Remark~\ref{rem:inftybeh} holds.
%\end{rem}

For the development and analysis of our algorithms, we use the following well-established and easy to verify result.
\begin{lem}\label{lem:VI_sol}
	Consider a mapping $\boldsymbol T(\cdot): \R^d \to \R^d$ and a convex closed set $Y \subseteq \R^d$. Given $\theta > 0$,
	\begin{align}
		\label{ex:projection_VI}
		\by^* \in \text{SOL}(Y, \boldsymbol T) \iff \by^* = \Proj_Y(\by^* - \theta \boldsymbol T(\by^*)).
	\end{align}
	%Given $\theta > 0$, if $\by^* = Proj_Y()$, then $\by^* \in SOL(Y,\boldsymbol T)$. \end{definition}
\end{lem}

\section{Proposed Optimization Algorithm}\label{sec:procedure}
\setlength{\textfloatsep}{15pt}
\begin{algorithm}[t!]
\caption{Payoff-based algorithm for learning  Nash equilibria in strongly monotone games}\label{alg:algorithm1}
\begin{algorithmic}[!t]
	\Require Action set $\Ab^i \subset \R^d$, the sequences $\{\sigma_t\}, \{\rho_t\}, \{\gamma_t\}$, initial state $\bmu^i(0)$.
	\For {$t = 0,1, \ldots$}
	\State  Sample $\bxi^i(t)$ according to probability density \eqref{eq:density}.
	\State Obtain the feasible action $\ab^i(t) = \Proj_{\Ab^i}[\bxi^i(t)]$.
	
	\noindent{\color{blue} { \footnotesize\ttfamily  /* Simultaneously and similarly,  other players choose their actions $\ab^{-i}(t)$.   /* }}
	\State Observe $ J^i(t) = J^i(\ab^1(t),\ldots,\ab^N(t))$.
	\State Perform the one-point gradient estimate $\boldsymbol m_1^i(t) = { J^i(t)}\frac{{\bxi^i(t)}  -\bmu^i(t)}{\sigma^2_t}$
	\State or perform the two-point gradient estimate $\boldsymbol m_2^i(t) = ( J^i(t) - J^i(\bmu^i(t)))\frac{{\bxi^i(t)} -\bmu^i(t)}{\sigma^2_t}$.
	\State Update the state: $\bmu^i(t+1)=\Proj_{(1-\rho_t)\Ab^i}[\bmu^i(t)-\gamma_t\boldsymbol m_j^i(t)]$, $j=1$ or $j=2$.
	
	\noindent{\color{blue} { \footnotesize\ttfamily  /* Simultaneously and similarly,  other players update their states $\bmu^{-i}(t+1)$.  /* }}
	\EndFor
	
\end{algorithmic}
\end{algorithm}
%We propose an algorithm for the online optimization problem that at each time step $t$ uses only the  experienced value of the function, $c_t(\bxi_t)$ for a played action $\bxi_t$. Then, we derive the regret bounds for our approach.

Let us denote by $\mb_j^i$, $j\in\{1,2\}$, some estimate of $\Mb^i$ in the pseudo-gradient of the game. The proposed method  to update the player $i$'s  so-called state $\bmu^i$ is as follows:
\begin{align}
\label{eq:alg}
\bmu^i(t+1)=\Proj_{(1-\rho_t)\Ab^i}[\bmu^i(t)-\gamma_t\mb_j^i(t)],
\end{align}
where $\bmu^i(0)\in \R^{Nd}$ is an arbitrary finite value, $\g_t$ is the step size,  and $\rho_t$ is a regularization parameter introduced to control the feasibility of the actions as will be detailed further. The step size $\gamma_t$ needs to be chosen based on the bias and variance of the pseudo-gradient estimates $\mb_j^i$, as well as the regularization parameters. The term $\mb_j^i(t)$, $j\in\{1,2\}$, is obtained using the payoff-based feedback as described below.

Given $\bmu^i(t)$, let  player $i$ sample the random vector $\bxi^i(t)$ according to the multidimensional normal distribution $\EuScript N(\bmu^i(t)=(\mu^i_1(t),\ldots,\mu^i_{d}(t))^{\top},\sigma_t)$ with the following density function:
\begin{align}\label{eq:density}
p^i&(\bx^i;\bmu^i(t),\sigma_{t})= \frac{1}{(\sqrt{2\pi}\sigma_{t})^{d}}\exp\left\{-\sum_{k=1}^{d}\frac{(x^i_k-\mu^i_k(t))^2}{2\sigma^2_{t}}\right\}.
\end{align}
Then, the action agent $i$ chooses is $\ab^i(t) = \Proj_{\Ab^i}[\bxi^i(t)]$.
In this work we consider two following set ups. The players can either have access only to the value of the cost function at the actual joint action (bandit setting), or they can use two estimations of their cost functions, namely at the currently played joint action $\ab(t)$ and also at the joint state $\bmu(t)$.
According to the bandit setting, the cost value $J^i(t)$ at the joint action $\ab(t)=(\ab^1(t),\ldots,\ab^N(t))\in \Ab$,  denoted by $J^i(t) =J^i(\ab(t)) $ is revealed to each agent $i$. Agent $i$ then estimates her local gradient $\frac{\partial J^i}{\partial \bmu^i}$ evaluated at the point of the joint state $\bmu(t)=(\bmu^1(t),\ldots,\bmu^N(t))$ as follows:
\begin{align}\label{eq:est_Gd}
\mb_1^i(t) = { J^i(t)}\frac{{\bxi^i(t)} -\bmu^i(t)}{\sigma^2_t}.
\end{align}
Otherwise, there is an extra peace of information available to each agent, namely $J^i(\bmu(t)) $. Then each agent uses the following estimation of the local gradient $\frac{\partial J^i}{\partial \bmu^i}$ at the point of $\bmu(t)$:
\begin{align}\label{eq:est_Gd2}
\mb_2^i(t) = (J^i(t) - J^i(\bmu(t)))\frac{{\bxi^i(t)} -\bmu^i(t)}{\sigma^2_t}.
\end{align}

The steps of the procedure run by each player are summarized in Algorithm \ref{alg:algorithm1}.
Before stating the convergence analysis,  let us provide insight into the procedure defined by Equation \eqref{eq:alg}  by deriving an analogy to a regularized stochastic gradient algorithm.

%{Characterizing the terms in the algorithm}
Let
\begin{align}\label{eq:densityfull}
p( \bx; \bmu, \sigma)=\prod_{i=1}^Np_i(x^i_1,\ldots,x^i_{d};\bmu^i,\sigma)
\end{align}
denote the  density function of the joint distribution of agents' states. Given $\sigma > 0$, for any $i\in[N]$ define $ \tilde{J}_i : \R^{Nd} \rightarrow \R$ as
\begin{align}
\label{eq:mixedJ}
\tilde{J}_i &(\bmu^1,\ldots,\bmu^N, \sigma)= \int_{\mathbb R^{Nd}}J_i(\bx)p( \bx; \bmu, \sigma)d\bx.
\end{align}
Thus, $\tilde{J}_i$, $i\in[N]$ is the $i$th player's cost function in mixed strategies sampled from the normal distribution with the density function in~\eqref{eq:densityfull}.
%We can now show that the second term inside the projection in \eqref{eq:regpl} is a sample of the gradient of this cost function $\tilde{J}_i$ with respect to the mixed strategies.
Let $\bmu(t) =(\bmu^1(t),\ldots,\bmu^N(t))$ and for $i\in[N]$, define $\tilde{\Mb}_i (\cdot)=(\tilde M_{i,1}(\cdot), \ldots, \tilde M_{i,d}(\cdot))^{\top}$
as the $d$-dimensional mapping with the following elements:
\begin{align}\label{eq:mapp2}
\tilde M_{i,k} (\bmu,\sigma)=\frac{\partial {\tilde J_i(\bmu, \sigma)}}{\partial \mu^i_k}, \mbox{ for $k\in[d]$}.
\end{align}

Using $\tilde{\Mb}_i (\cdot)$, $i\in[N]$ defined above
we can rewrite the algorithm step in \eqref{eq:alg} in the form below:
\begin{align}
\label{eq:pbavmu}
\bmu^i(t+1) =\Proj_{(1-\rho_t)A_i}[&\bmu^i(t) -\gamma_t\sigma^2_t\big(\tilde\Mb_i(\bmu(t)) +\Rb_{i,j}(\bxi(t),\bmu(t),\sigma_t)+\Pb_i(\bxi(t),\bmu(t),\sigma_t)\big)],
\end{align}
for all $i\in[N]$, where  $\Rb_{i,j}$, $j=1,2$, and $\Pb_i$ are defined as follows:
\begin{align}
&\Pb_i(\bxi(t),\bmu(t),\sigma_t) = \frac{\bxi^i(t) -\bmu^i(t)}{\sigma^2_t}(J_i(\ba(t))-J_i(\bxi(t)))\label{eq:Pterm},\\
&\Rb_{i,j}(\bxi(t),\bmu(t),\sigma_t) = \tilde{\mb}_j^i(t) - \tilde{\Mb}_i (\bmu(t),\sigma_t), \, \mbox{where}\label{eq:Rterm} \cr
&\tilde{\mb}_j^i(t) =\begin{cases}
	J_i(\bxi(t))\frac{\bxi^i(t) -\bmu^i(t)}{\sigma^2_t}, \, &\mbox{if $j=1$ (one-point estimation)},\\
	(J_i(\bxi(t)) - J_i(\bmu(t)))\frac{\bxi^i(t) -\bmu^i(t)}{\sigma^2_t}, \, &\mbox{if $j=2$ (two-point estimation)}.
\end{cases}
\end{align}

The following lemmas demonstrate the behavior of the terms $\Pb_i$ and $\Rb_{i,j}$. Due to the sampling from the Gaussian distribution with unbounded support, we need the following assumption on the cost functions' behavior at infinity.
\begin{assumption}	\label{assum:infty}
	Each function $J^i(\bx) = O(\exp\{\|\bx\|^{\alpha}\})$ as $\|\bx\|\to\infty$, where $\alpha<2$.
\end{assumption}

First, we demonstrate that the mapping $\tilde{\Mb}_i (\bmu(t))$ evaluated at $\bmu(t)$ is equivalent to the expected pseudo-gradient. That is,
\begin{align}\label{eq:gradmix}
\tilde{\Mb}_i (&\bmu(t))=\int_{\mathbb R^{Nd}}{\Mb_i} (\bx)p(\bx;\bmu(t),\sigma_t)d\bx.
\end{align}
Moreover, this lemma proves that $\tilde{\Mb}_i (\bmu(t))$ is equal to the expectation of the term \[\tilde{\mb}_j^i(t) = (\tilde{m}_j^{i,1}(t), \ldots, \tilde{m}_j^{i,d}(t))\in\R^d, \quad j = 1, 2.\]

\begin{lem}\label{lem:sample_grad}
Given Assumptions\r\ref{assum:convex} and~\ref{assum:infty}, for $j=1,2$,
\begin{align}\label{eq:deriv}
	&\tilde M_{i,k} (\bmu(t)) =\E\{\tilde{m}_j^{i,k}(t)|\xi^i_k(t)\sim\EuScript N(\mu_k^i(t),\sigma_t), i\in[N], k\in[d]\} \cr
	&\qquad =\E\{M_{i,k}(\bxi^1,\ldots,\bxi^N)|\xi^i_k(t)\sim\EuScript N(\mu_k^i(t),\sigma_t), i\in[N], k\in[d]\}.
\end{align}
\end{lem}
Proof is in Appendix~\ref{app:sample}.

Lemma~\ref{lem:sample_grad} above implies that
\begin{align}
\label{eq:mathexp2}
\Rb_{i,j}(\bxi(t),\bmu(t),\sigma_t) = \tilde{\mb}_j^i(t) - \E_{\bxi(t)}\{\tilde{\mb}_j^i(t)\},\; \;i\in[N].
\end{align}

Further, for the sake of notation simplicity, we may use $\Rb(t) = \Rb(\bxi(t), \bmu(t),\sigma_t)$.
Let $\EuScript F_t$ be the $\sigma$-algebra generated by the random variables $\{\bmu(k),\bxi(k)\}_{k\le t}$.
For the second moment of the term  $\Rb$ we prove the following lemma.
\begin{lem}\label{lem:Rsq}
Under Under Assumptions~\ref{assum:convex} and~\ref{assum:infty}, as $\s_t\to 0$,
\[\E\{\|\Rb_{i,j}(t)\|^2 | \EuScript F_t\} = \begin{cases}
	O\left(\frac{1}{\sigma_t^2}\right), \, &\mbox{ if $j=1$},\\
	O(1). \, &\mbox{ if $j=2$}.
\end{cases}\]
\end{lem}
Proof is in Appendix~\ref{app:Rsq}.

Finally, the term
\begin{align*}
\Pb(\bxi(t), \bmu(t),\sigma_t) = (\Pb_1(\bxi(t),\bmu(t),\sigma_t),\ldots,\Pb_N(\bxi(t), \bmu(t),\sigma_t)),
\end{align*}
is the vector of the difference between the gradient estimation based on the state $\bxi(t)\in \R^{Nd}$ and the  played action $\ba(t)\in\Ab$.
To characterize this term we prove the following statement.
\begin{lem}\label{lem:projTerm}
Under Assumptions~\ref{assum:convex} and~\ref{assum:infty} we have
\begin{align}
	\E\left\{\|\Pb_i(\bxi(t), \bmu(t),\sigma_t)\|\;|\EuScript F_t\right\}&= \E\left\{|J_i(\ab(t))-J_i(\bxi(t))|\frac{\|\bxi^i(t) -\bmu^i(t)\|}{\sigma^2_t}\;|\EuScript F_t\right\}\cr
	&=O\left(\left(\frac{e^{-\frac{{\rho}_t^2}{2\sigma_t^2}}}{\sigma_t^{Nd+1}}\right)^{\frac{1}{2}}\right) \mbox{ almost surely. }
\end{align}
\end{lem}
Proof is in Appendix~\ref{app:proj}.

%\begin{rem}\label{rem:Lemmas}
%According to Remark~\ref{rem:AssumImpl}, Lemmas~\ref{lem:sample_grad}-\ref{lem:projTerm} above hold if Assumption\r\ref{assum:LipschitzTwoP} and~\ref{assum:infty} therein are replaced by Assumption~\ref{assum:Lipschitz}.
%	\end{rem}

\section{Convergence Analysis}\label{sec:rate}
\subsection{One-point Estimations ($j=1$)}
We will provide the analysis of Algorithm~\ref{alg:algorithm1} in the case $j=1$ (one-point gradient estimations) under the following smoothness assumption.
\begin{assumption}	\label{assum:Lipschitz}
	The pseudo-gradient $\Mb:\R^{Nd}\to\R^{Nd}$, defined in Assumption~\ref{assum:CG_grad} fulfills one of the following conditions:
	\begin{enumerate}
		\item $\Mb$ is twice  differentiable over $\R^{Nd}$; \label{num:11}
		\item $\Mb$ is differentiable and Lipschitz continuous over $\R^{Nd}$ with some constant $K$. \label{num:12}
	\end{enumerate}
\end{assumption}

\begin{rem}
	 Note that under Assumption~\ref{assum:Lipschitz}.\ref{num:12} each function $J_i(\bx) = O(\|\bx\|^{2})$ as $\|\bx\|\to\infty$ and, thus, Assumption~\ref{assum:infty} holds in this case.
\end{rem}

In the case of one-point gradient estimations, we will base our analysis on the algorithm's representation in~\eqref{eq:pbavmu}. Thus, in this subsection we exploit the properties of the term $\tilde{\Mb}_i (\bmu(t))$ therein.

Let us now focus on the mapping $\tilde{\Mb}^{(t)}(\cdot) =(\tilde{\Mb}^{(t)}_1 (\cdot), \ldots, \tilde{\Mb}^{(t)}_N (\cdot))$, where for any $\bmu\in\R^{Nd}$
\begin{align}
\label{eq:smooth_t_game}
\tilde{\Mb}^{(t)}_i (\bmu)=\int_{\mathbb R^{Nd}}{\Mb_i} (\bx)p(\bx;\bmu,\sigma_t)d\bx
\end{align}
for some given $\sigma_t$.
Note that if  $\bmu=\bmu(t)$, where $\bmu(t)$ is updated according to~\eqref{eq:alg}, then $\tilde{\Mb}^{(t)}_i (\bmu)$ is equal to $\tilde{\Mb}_i (\bmu(t))$, according to Lemma~\ref{lem:sample_grad} (see also \eqref{eq:gradmix}). Such mapping is the pseudo-gradient in the mixed strategies, given that the joint action is generated by the normal distribution with the density~\eqref{eq:densityfull}.

Our next results describe the properties of the mapping $\tilde{\Mb}^{(t)}$ which are important for the further algorithm's  analysis.

First, we note that the following statement, which describes the connection between ${\Mb}^{(t)}$ and ${\Mb}^{(t-1)}$ as well as ${\Mb}^{(t)}$ and $\Mb$, takes place.
\begin{prop}\label{prop:tildeM}
Let the pseudo-gradient $\Mb$ be differentiable over $\R^{Nd}$. Then, under Assumptions~\ref{assum:infty}, for any $\bmu\in\Ab$, we have
\begin{align}\label{eq:tildeM}
	\|\tilde{\Mb}^{(t)} (\bmu)-\tilde{\Mb}^{(t-1)}(\bmu)\| = O\left(|\sigma_t - \sigma_{t-1}|\right),
\end{align}
\begin{align}\label{eq:MvsTildeM}
	\|\tilde{\Mb}^{(t)} (\bmu)-\Mb(\bmu)\| = O(\sigma_t).
\end{align}
\end{prop}
See Appendix~\ref{app:convex} for the proof.

\begin{prop}\label{prop:strMon}
Let Assumptions~\ref{assum:CG_grad}--\ref{assum:Lipschitz} hold. Then  $\tilde M (\bmu)$ is continuous and differentiable over $\R^{Nd}$ and Lipschitz continuous over any compact set. Moreover, let $\sigma_t\to 0$ as $t\to\infty$.  Then there exists some finite $T$ such that for each $t\ge T$ the mapping $\tilde{\Mb}^{(t)}$ is strongly monotone with the constant $\frac{\nu}{2}$.
%and any $\bmu_1$, $\bmu_2\in\Ab$
%\[(\tilde{\Mb}^{(t)}(\bmu_1)-\tilde{\Mb}^{(t)}(\bmu_2),\bmu_1-\bmu_2)\ge \frac{\nu}{2}\|\bmu_1-\bmu_2\|^2.\]
\end{prop}
See Appendix~\ref{app:strMon} for the proof.

Proposition above states that under Assumptions~\ref{assum:CG_grad}--\ref{assum:Lipschitz} the mapping $\tilde{\Mb}^{(t)}$ is continuous over $\R^{Nd}$. Thus, if Assumption~\ref{assum:compact} holds,  the set $SOL(\Ab, \tilde{\Mb}^{(t)})$ is not empty (see~\cite{FaccPang1}). The following proposition provides the connection between Nash equilibria in the game $\Gamma(N,\{A_i\},\{J_i\})$ and the set $SOL(\Ab, \tilde{\Mb}^{(t)})$.

\begin{prop}\label{prop:distSol}
Let Assumptions~\ref{assum:convex}--\ref{assum:infty} hold. Let $\ab^*$ be the unique Nash equilibrium in $\Gamma(N,\{A_i\},\{J_i\})$ and $\bmu^*_t\in SOL(\Ab, \tilde{\Mb}^{(t)})$. Then, as $\sigma_t\to 0$,
\[\|\bmu^*_t-\ab^*\| = O(\sigma_t).\]
\end{prop}
See Appendix~\ref{app:distSol} for the proof.

The next proposition we need to further analyze convergence of the proposed algorithm relates two points $\bmu^*_t$ and $\bmu^*_{t-1}$ from the sets $SOL(\Ab, \tilde{\Mb}^{(t)})$ and $SOL(\Ab, \tilde{\Mb}^{(t-1)})$ respectively.

\begin{prop}\label{prop:tVSt-1}
Let Assumptions~\ref{assum:CG_grad}--\ref{assum:Lipschitz} hold and $\sigma_t\to 0$ as $t\to\infty$.
Then, as $t\to\infty$,
\[\|\bmu^*_t-\bmu^*_{t-1}\|= O(|\sigma_t-\sigma_{t-1}|).\]
\end{prop}
See Appendix~\ref{app:tVSt-1} for the proof.

Let us also consider  the following auxiliary sequence $\yb(t)$:
\begin{align}\label{eq:y}
%\zb(t+1) & = \Prt[\zb(t)-\gamma_t\Mb(\zb(t))], \\
\yb(t) & = SOL((1-\rho_t)\Ab, \tilde{\Mb}^{(t)}),
\end{align}
where $\rho_t\in[0,1)$.
We prove the following statements regarding the sequence $\yb(t)$.

\begin{prop}\label{prop:diff}
Let Assumptions~\ref{assum:CG_grad}--\ref{assum:Lipschitz} hold and $\sigma_t\to 0$ as $t\to\infty$.
Then, as $t\to\infty$,
\[\|\yb(t)-\yb(t-1)\|=  O(|\sigma_t-\sigma_{t-1}|)+ O(|\rho_t-\rho_{t-1}|).\]
\end{prop}
See Appendix~\ref{app:diff} for the proof.

\begin{prop}\label{prop:y-mu}
Let Assumptions~\ref{assum:CG_grad}-~\ref{assum:Lipschitz} hold and $\sigma_t\to 0$ as $t\to\infty$.
 Then, as $t\to\infty$,
\[\|\yb(t)-\bmu_t^*\| =  O(\rho_t).\]
\end{prop}
See Appendix~\ref{app:y-mu} for the proof.

With all these results at place, we are ready to formulate the main result.
\begin{theorem}\label{th:main}
	Let the states $\bmu^i(t)$, $i\in[N]$, evolve according to Algorithm~\ref{alg:algorithm1} with the gradient estimations $\mb_1^i(t)$.
	Let Assumptions\r\ref{assum:convex}--\ref{assum:Lipschitz} hold.
	Moreover, let the parameters in the procedure be chosen as follows: $\gamma_t=\frac{4}{\nu t}$, $\sigma_t = \frac{a}{t^{\frac{1}{4}}}$, $\rho_t = \frac{1}{t^{\frac{1}{4}-\varepsilon}}$ for arbitrary small $\varepsilon>0$. Here $\nu$ is the strong monotonicity constant from Assumption~\ref{assum:CG_grad}.
	
	Then the joint state $\bmu(t)$ converges almost surely to the unique Nash equilibrium $\bmu^*=\ba^*$ of the game $\Gamma$, whereas the joint state $\bxi(t)$ converges in probability to $\ba^*$. Moreover,
	\begin{align*}
		\E \|\bmu(t) - \ba^*\|^2=O\left(\frac{1}{t^{1/2-\varepsilon}}\right).
	\end{align*}
\end{theorem}

\begin{proof}
	Let us notice that due to the theorem's conditions and the particular choice $\sigma_t = \frac{a}{t^{\frac{1}{4}}}\to 0$, as $t\to\infty$, Propositions~\ref{prop:strMon}-\ref{prop:y-mu} hold.
	
	We consider $\|\bmu(t+1)-\by(t)\|^2$, where the sequence $\by(t)$ is defined in~\eqref{eq:y}.
	We aim to bound the growth of $\|\bmu(t+1)-\by(t)\|^2$ in terms of $\|\bmu(t)-\by(t-1)\|^2$ and, thus, to obtain the convergence rate of the sequence  $\|\bmu(t+1)-\by(t)\|^2$. Further, we aim to apply Propositions~\ref{prop:y-mu} and~\ref{prop:distSol} estimating the distances $\|\by(t-1)-\bmu^*_t\|$  and $\|\bmu^*_t-\ab^*\|$ respectively to conclude the result.
	
	%For this purpose let $\e(t)$ be some positive sequence such that
	%\begin{align}\label{eq:e}
	%\sum_{t=1}^{\infty}\e(t)<\infty \, \mbox{ and } \, \sum_{t=1}^{\infty}\left(1+\frac{1}{\e(t)}\right)(\rho_t-\rho_{t-1})^2<\infty.
	%\end{align}
	%Under Assumption~\ref{assum:timestep} the choice of such $\e(t)$ is possible. Indeed, as $(\rho_t-\rho_{t-1})^2\sim \frac{1}{t^{2r+2}}$, the sequence $\e(t)$ can be set up as $\e(t) = \frac{1}{t^{1+e}}$, where $e<2r$.
	%Under such choice of $\e(t)$ we have $\e(t)\rightarrow 0$ and $\left(1+\frac{1}{\e(t)}\right)(\rho_t-\rho_{t-1})^2 \rightarrow 0$ as $t\to\infty$. Hence,  $\forall \bmu\in\Ab$
	%   \begin{align}\label{eq:connection1}
		%    \|\bmu-\by(t)\|^2 \le& (1+\e(t))\|\bmu - \by(t-1)\|^2  +\left(1 + \frac{1}{\e(t)}\right)\|\by(t-1)-\by(t)\|^2 \\ \nonumber
		% \le& (1+\e(t))\|\bmu - \by(t-1)\|^2  \\ \nonumber
		% & +\left(1 + \frac{1}{\e(t)}\right)O\left(\frac{1}{t^2}+(\rho_{t-1}-\rho_t)^2\right),
		%   \end{align}
	%where in the last inequality we used Lemma~\ref{prop:diff}.
	We analyze each term in the following sum $\sum_{i=1}^{N} \|\bmu^i(t+1)-\by^i(t)\|^2$.
	From the procedure for the update of $\bmu(t)$ in ~\eqref{eq:pbavmu}, the non-expansion property of the projection operator, the fact that $\by(t)$ belongs to $SOL((1-\rho_t)\Ab,\tilde{\Mb}^{(t)}(\by))$, namely, that $ \forall i\in[N]$
	\[\by^i(t) = \Proj_{(1-\rho_t)A_i}[\by^i(t) - \gamma_t(\tilde{\Mb}^{(t)}_i(\by(t))],\]
	where $\tilde{\Mb}^{(t)}=\tilde{\Mb}$,
	we obtain that for any $i\in[N]$
	\begin{align}\label{eq:nonexp}
		\|&\bmu^i(t+1)-\by^i(t)\|^2 \le \|\bmu^i(t)-\by^i(t) -\gamma_t\big[\tilde{\Mb}_i(\bmu(t)) - \tilde{\Mb}_i(\by(t))\cr
		& +\Rb_{i,j}(\bxi(t),\bmu(t))+\Pb_i(\bxi(t),\bmu(t))\big]\|^2\cr
		& = \|\bmu^i(t)-\by^i(t)\|^2 \cr
		&\qquad- 2\gamma_t(\tilde{\Mb}_i(\bmu(t))- \tilde{\Mb}_i(\by(t)), \bmu^i(t)-\by^i(t)) \cr
		&\qquad-2\gamma_t(\Rb_{i,j}(t), \bmu^i(t)-\by^i(t)) \cr
		&\qquad-2\gamma_t(\Pb_i(t), \bmu^i(t)-\by^i(t)) \cr
		&\qquad + \gamma^2_t\|\Gb_i(t)\|^2,
	\end{align}
	where, for ease of notation, we have defined $\Rb(t) = \Rb(\bxi(t), \bmu(t),\sigma_t)$, $\Pb(t) = \Pb(\bxi(t), \bmu(t),\sigma_t)$, and
	\begin{align}\label{eq:G_1}
		\Gb_i(t) = &\tilde{\Mb}_i(\bmu(t)) - \tilde{\Mb}_i(\by(t)) \cr
		&+\Rb_{i,j}(t)+\Pb_i(t).
	\end{align}
	We expand $\Gb_i$ as below and bound the terms in the expansion.
	\begin{align}\label{eq:G_1_norm}
		&\|\Gb_i(t)\|^2 = \|\tilde{\Mb}_i(\bmu(t)) - \tilde{\Mb}_i(\by(t))\|^2 \cr
		&+ \|\Rb_{i,j}(t)\|^2+\|\Pb_i(t)\|^2\cr
		& + 2(\Pb_i(t),\Rb_{i,j}(t))\cr
		& +2(\tilde{\Mb}_i(\bmu(t)) - \tilde{\Mb}_i(\by(t)),\Rb_{i,j}(t)+\Pb_i(t)).
	\end{align}
	Thus, by taking into account~\eqref{eq:mathexp2}, which implies $\E\{\Rb_{i,j}(t)|\EuScript F_t\} =\boldsymbol 0$ for any $\bmu$, and the Cauchy-Schwarz inequality, we get from \eqref{eq:nonexp}
	\begin{align}\label{eq:firstEst}
		&\E\{\|\bmu^i(t+1)-\by^i(t)\|^2|\EuScript F_t\} \cr
		&\le \|\bmu^i(t)-\by^i(t)\|^2 \cr
		&\quad- 2\gamma_t(\tilde{\Mb}_i(\bmu(t))- \tilde{\Mb}_i(\by(t)), \bmu^i(t)-\by^i(t)) \cr
		&\quad + 2\gamma_t\E\{\|\Pb_i(t)\||\EuScript F_t\}\|\bmu^i(t)-\by^i(t)\|\cr
		&\quad + \gamma^2(t)\E\{\|\Gb_i(t)\|^2|\EuScript F_t\}\cr
		& \le\|\bmu^i(t)-\by^i(t)\|^2 \cr
		&\quad- 2\gamma_t(\tilde{\Mb}_i(\bmu(t))- \tilde{\Mb}_i(\by(t)), \bmu^i(t)-\by^i(t)) \cr
		&\quad + 2\gamma_t\E\{\|\Pb_i(t)\||\EuScript F_t\}\|\bmu^i(t)-\by^i(t)\|\cr
		&\quad + \gamma^2(t)[\|\tilde{\Mb}_i(\bmu(t)) - \tilde{\Mb}_i(\by(t))\|^2\cr
		&\quad+\E\{\|\Rb_{i,j}(t)\|^2+\|\Pb_i(t)\|^2|\EuScript F_t\}\cr
		&\quad+2\E\{ \|\Pb_i(t)\|\|\Rb_{i,j}(t)\||\EuScript F_t\}\cr
		&\quad +2(\|\tilde{\Mb}_i(\bmu(t)) - \tilde{\Mb}_i(\by(t))\|)\E\{\|\Pb_i(t)\||\EuScript F_t\}].
	\end{align}
	Next, taking into account Lemmas~\ref{lem:Rsq},~\ref{lem:projTerm} and the fact that $\frac{\s_t}{\rho_t}\to 0$ as $t\to\infty$ we conclude that $\E\{\|\Pb_i(t)\||\EuScript F_t\}$ and $\E\{\|\Pb_i(t)\|^2|\EuScript F_t\}$ decrease geometrically fast in comparison with other ones. Hence, we obtain  from~\eqref{eq:firstEst} that
	\begin{align}
		&\E\{\|\bmu^i(t+1)-\by^i(t)\|^2|\EuScript F_t\} \cr
		&\le\|\bmu^i(t)-\by^i(t)\|^2 \cr
		&\quad- 2\gamma_t(\tilde{\Mb}_i(\bmu(t))- \tilde{\Mb}_i(\by(t)), \bmu^i(t)-\by^i(t))\cr
		& \quad+ \gamma^2(t)[\|\tilde{\Mb}_i(\bmu(t)) - \tilde{\Mb}_i(\by(t))\|^2+\E\{\|\Rb_{i,j}(t)\|^2|\EuScript F_t\}].
	\end{align}
	Thus, using continuity of $\tilde{\Mb}$ (Proposition~\ref{prop:strMon}) and, thus, boundedness of $\|\tilde{\Mb}_i(\bmu(t)) - \tilde{\Mb}_i(\by(t))\|$ as well as taking into account the result  of  Lemma~\ref{lem:Rsq}, we conclude that
	\begin{align}\label{eq:secondEst}
		&\E\{\|\bmu^i(t+1)-\by^i(t)\|^2|\EuScript F_t\} \cr
		&\le\|\bmu^i(t)-\by^i(t)\|^2 \cr
		&\quad- 2\gamma_t(\tilde{\Mb}_i(\bmu(t))- \tilde{\Mb}_i(\by(t)), \bmu^i(t)-\by^i(t))+ h(t),
	\end{align}
	where
	\begin{align}\label{eq:ht}
		h(t) = 	O\left(\frac{\gamma^2_t}{\sigma^2_t}\right).
	\end{align}
	Thus, due to Proposition~\ref{prop:strMon}, we conclude from~\eqref{eq:secondEst} by summing up the inequalities over $i=1,\ldots,N$,
	\[\E\{\|\bmu(t+1)-\by(t)\|^2|\EuScript F_t\}\le(1-\nu\gamma_t)\|\bmu(t)-\by(t)\|^2+N h(t).\]
	Next, by using the inequality
	\[\|\bmu(t)-\by(t)\|^2\le (1+\epsilon_t)\|\bmu(t)-\by(t-1)\|^2 + (1+1/\epsilon_t)\|\by(t-1) -\by(t)\|^2,\]
	which holds for any $\epsilon_t>0$,
	we get
	\begin{align*}
		\E\{\|\bmu(t+1)-\by(t)\|^2|\EuScript F_t\}&\le(1-\nu\gamma_t)(1+\epsilon_t)\|\bmu(t)-\by(t-1)\|^2 \cr
		&\,  + (1-\nu\gamma_t)(1+1/\epsilon_t)\|\by(t-1) -\by(t)\|^2+h(t)\cr
		&\le\left(1-\frac{\nu\gamma_t}{2}\right)\|\bmu(t)-\by(t-1)\|^2 \cr
		&\qquad+ O\left(\frac{|\sigma_t-\sigma_{t-1}|^2+|\rho_t-\rho_{t-1}|^2}{\gamma_t}+h(t)\right),
	\end{align*}
	where in the last inequality we chose $\epsilon_t = \frac{\nu\gamma_t}{2}$ and used Proposition~\ref{prop:diff}.
	Thus, given the settings for the parameter $\g_t = \frac{4}{\nu t}$, the definition of $h(t)$ (see~\eqref{eq:ht}),  we conclude that
	\begin{align}\label{eq:LV2}
		\E\{\|\bmu(t+1)-\by(t)\|^2|\EuScript F_t\}\le\left(1-\frac{2}{t}\right)\|\bmu(t)-\by(t-1)\|^2 + O\left(\frac{1}{t^{3/2-\varepsilon}}\right),
	\end{align}
	Here we used the definition of $h_t$ from~\eqref{eq:ht}, which implies $h_t = O\left(\frac{1}{t^{3/2}}\right)$. Moreover,
	\begin{align*}
		\frac{|\sigma_t-\sigma_{t-1}|^2+|\rho_t-\rho_{t-1}|^2}{\gamma_t} = 		O\left(\frac{1}{t^{3/2-\varepsilon}}\right), \quad &\mbox{if $j=1$}.
	\end{align*}
	Thus, \eqref{eq:LV2} implies that $\bmu(t)$ converges to $\by(t-1)$ almost surely (see Chung's Lemma 4 in Chapter 2.2~\cite{Polyak}, for the reader's convenience we provide this lemma in Appendix~\ref{app:Ch}). Moreover, due to Propositions~\ref{prop:y-mu} and~\ref{prop:distSol}, the almost sure convergence of $\bmu(t)$ to $\ab^*$ takes place. Taking into account that $\bxi(t)\sim\EuScript N(\bmu(t),\sigma_t)$ and $\sigma_t\to 0$ as $t\to\infty$, we conclude that $\bxi(t)$ converges weakly to a Nash equilibrium $\ab^*$. Moreover, according to the Portmanteau Lemma \cite{portlem}, this convergence is also in probability.
	Next, by taking the full expectation of the both sides in~\eqref{eq:LV2}, we obtain
	\begin{align*}
		&\E[\|\bmu(t+1)-\by(t)\|^2]\le\left(1-\frac{2}{t}\right)\E\|\bmu(t) - \by(t-1)\|^2+ O\left(\frac{1}{t^{p}}\right).
	\end{align*}
	By applying the Chung's lemma (see Appendix~\ref{app:Ch}) to the inequality above, we conclude that
	\[\E[\|\bmu(t+1)-\by(t)\|^2] = O\left(\frac{1}{t^{1/2-\varepsilon}}\right).\]
	Finally,
	\begin{align}\label{eq:final}
		\E[\|\bmu(t+1)-\ab^*\|^2]  & \le  2\E[\|\bmu(t+1)-\by(t)\|^2] + 2\|\by(t)-\bmu^*_t\|^2  + 2\|\bmu^*_t-\ab^*\|^2.
	\end{align}
	Thus, applying Propositions~\ref{prop:distSol} and~\ref{prop:y-mu} with the given choice of the parameters, we obtain the result.
\end{proof}

%%%%%%%%%%%%%%%%%%%%%%%%%%%%%%%5
%%%%%%%%%%%%%%%%%%%%%%%%%%%%%
%%%TWO%%%%%%%%%%%%%%%%%%%%%%%%%
%%%%%%%%%%%%%%%%%%%%%%%%%%%
%%%%%%%%%%%%%%%%%%%%%%%%%%%%%%

\subsection{Two-point Estimations ($j=2$)}
We will provide the analysis of Algorithm~\ref{alg:algorithm1} in the case $j=2$ (two-point gradient estimations) under the following smoothness assumption.
\begin{assumption}	\label{assum:Lipschitz2}
	The pseudo-gradient $\Mb:\R^{Nd}\to\R^{Nd}$, defined in Assumption~\ref{assum:CG_grad} fulfills one of the following conditions:
	\begin{enumerate}
		\item $\Mb$ is differentiable over $\R^{Nd}$; \label{num:21}
		\item $\Mb$ is Lipschitz continuous over $\R^{Nd}$ with some constant $l$. \label{num:22}
	\end{enumerate}
\end{assumption}

For the case, when two points are used by agents for gradient estimations ($j=2$ in Algorithm~\ref{alg:algorithm1}), we focus on the following representation of the proposed learning procedure from~\eqref{eq:pbavmu}:
\begin{align}
	\label{eq:pbavmu1}
	\bmu^i(t+1) =\Proj_{(1-\rho_t)A_i}[&\bmu^i(t) -\gamma_t\sigma^2_t\big(\Mb_i(\bmu(t)) +\cr
	&+\Qb_i(\bmu(t)) +\Rb_{i,2}(\bxi(t),\bmu(t),\sigma_t)+\Pb_i(\bxi(t),\bmu(t),\sigma_t)\big)],
\end{align}
where $\Qb_i(\bmu(t)) = \tilde\Mb_i(\bmu(t)) - \Mb_i(\bmu(t))$.
That is why to analyze the procedure in this case by comparing its iterates with ones based on the actual gradients, we introduce the following sequence:
\begin{align}\label{eq:z}
	\zb(t+1) & = \Prt[\zb(t)-\gamma_t\Mb(\zb(t))],
\end{align}
where $\rho_t\in[0,1)$.
For the forthcoming analysis we will need the following statements regarding the sequence $\zb(t)$.

\begin{prop}\label{prop:diffz}
	Let Assumptions~\ref{assum:CG_grad} and \ref{assum:Lipschitz2} hold.
	Then
	\[\|\zb(t)-\zb(t-1)\|= O(|\rho_t-\rho_{t-1}|).\]
\end{prop}
%See Appendix~\ref{app:diffz} for the proof.

\begin{prop}\label{prop:z-mu}
	Let Assumptions~\ref{assum:CG_grad} and \ref{assum:Lipschitz2} hold. Let $\ab^*$ be the unique Nash equilibrium in $\Gamma(N,\{A_i\},\{J_i\})$.
	 Then
	\[\|\zb(t)-\ab^*\| =  O(\rho_t).\]
\end{prop}
Note that the proofs of Propositions~\ref{prop:diffz},~\ref{prop:z-mu} are analogous to ones of Propositions~\ref{prop:diff},~\ref{prop:y-mu} respectively, where the mapping $\tilde\Mb^{(t)}$ is replaced by the pseudo-gradient $\Mb$ with the corresponding properties.

The following theorem formulates the main result on the convergence and its rate.
\begin{theorem}\label{th:main2}
	Let the states $\bmu^i(t)$, $i\in[N]$, evolve according to Algorithm~\ref{alg:algorithm1} with the gradient estimations $\mb_j^i(t)$, $j=2$.
	Let Assumptions\r\ref{assum:convex}--\ref{assum:infty}, and~\ref{assum:Lipschitz2} hold.
	Moreover, let the parameters in the procedure be chosen as follows: $\gamma_t=\frac{4}{\nu t}$, $\sigma_t = \frac{b}{t^s}$, $\rho_t = \frac{c}{t^r}$, $1\le r<s$. Here $\nu$ is the strong monotonicity constant from Assumption~\ref{assum:CG_grad}.
	
	Then the joint state $\bmu(t)$ converges almost surely to the unique Nash equilibrium $\bmu^*=\ba^*$ of the game $\Gamma$, whereas the joint state $\bxi(t)$ converges in probability to $\ba^*$. Moreover,
	\begin{align*}
		\E \|\bmu(t) - \ba^*\|^2=O\left(\frac{1}{t}\right).
	\end{align*}
\end{theorem}

\begin{proof}
	We follow here the logic of the proof for Theorem~\ref{th:main}. However, in this case we consider $\|\bmu(t+1)-\zb(t)\|^2$, where the sequence $\zb(t)$ is defined in~\eqref{eq:z}.
	
	%For this purpose let $\e(t)$ be some positive sequence such that
	%\begin{align}\label{eq:e}
	%\sum_{t=1}^{\infty}\e(t)<\infty \, \mbox{ and } \, \sum_{t=1}^{\infty}\left(1+\frac{1}{\e(t)}\right)(\rho_t-\rho_{t-1})^2<\infty.
	%\end{align}
	%Under Assumption~\ref{assum:timestep} the choice of such $\e(t)$ is possible. Indeed, as $(\rho_t-\rho_{t-1})^2\sim \frac{1}{t^{2r+2}}$, the sequence $\e(t)$ can be set up as $\e(t) = \frac{1}{t^{1+e}}$, where $e<2r$.
	%Under such choice of $\e(t)$ we have $\e(t)\rightarrow 0$ and $\left(1+\frac{1}{\e(t)}\right)(\rho_t-\rho_{t-1})^2 \rightarrow 0$ as $t\to\infty$. Hence,  $\forall \bmu\in\Ab$
	%   \begin{align}\label{eq:connection1}
		%    \|\bmu-\by(t)\|^2 \le& (1+\e(t))\|\bmu - \by(t-1)\|^2  +\left(1 + \frac{1}{\e(t)}\right)\|\by(t-1)-\by(t)\|^2 \\ \nonumber
		% \le& (1+\e(t))\|\bmu - \by(t-1)\|^2  \\ \nonumber
		% & +\left(1 + \frac{1}{\e(t)}\right)O\left(\frac{1}{t^2}+(\rho_{t-1}-\rho_t)^2\right),
		%   \end{align}
	%where in the last inequality we used Lemma~\ref{prop:diff}.
	We analyze each term in the following sum $\sum_{i=1}^{N} \|\bmu^i(t+1)-\zb^i(t)\|^2$.
	From the procedure for the update of $\bmu(t)$ in ~\eqref{eq:pbavmu1}, the non-expansion property of the projection operator, the fact that $\zb(t)$ belongs to $SOL((1-\rho_t)\Ab,{\Mb}(\zb))$, namely, that $ \forall i\in[N]$
	\[\zb^i(t) = \Proj_{(1-\rho_t)A_i}[\zb^i(t) - \gamma_t({\Mb}(\zb))],\]
	we obtain that for any $i\in[N]$
	\begin{align}\label{eq:nonexp}
		\|&\bmu^i(t+1)-\zb^i(t)\|^2 \le \|\bmu^i(t)-\zb^i(t) -\gamma_t\big[{\Mb}_i(\bmu(t)) - {\Mb}_i(\zb(t))\cr
		& \Qb_i(\bmu(t))+\Rb_{i,2}(\bxi(t),\bmu(t))+\Pb_i(\bxi(t),\bmu(t))\big]\|^2\cr
		& = \|\bmu^i(t)-\zb^i(t)\|^2   -2\gamma_t(\Qb_{i}(t), \bmu^i(t)-\zb^i(t))\cr
		&\qquad- 2\gamma_t({\Mb}_i(\bmu(t))- {\Mb}_i(\zb(t)), \bmu^i(t)-\zb^i(t)) \cr
		&\qquad-2\gamma_t(\Rb_{i,2}(t), \bmu^i(t)-\zb^i(t)) \cr
		&\qquad-2\gamma_t(\Pb_i(t), \bmu^i(t)-\zb^i(t)) \cr
		&\qquad + \gamma^2_t\|\Hb_i(t)\|^2,
	\end{align}
	where, for ease of notation, we have defined $\Qb_{i}(t) = \Qb_{i}(\bmu(t))$, $\Rb(t) = \Rb(\bxi(t), \bmu(t),\sigma_t)$, $\Pb(t) = \Pb(\bxi(t), \bmu(t),\sigma_t)$, and
	\begin{align}\label{eq:G_1}
		\Hb_i(t) = &\tilde{\Mb}_i(\bmu(t)) - \tilde{\Mb}_i(\by(t)) \cr
		& +\Qb_i(t)+\Rb_{i,j}(t)+\Pb_i(t).
	\end{align}
 Taking into account that  $\Qb_i(\bmu(t)) = \tilde\Mb_i(\bmu(t)) - \Mb_i(\bmu(t))$, $\tilde\Mb_i(\bmu(t)) = \tilde\Mb^{(t)}_i(\bmu(t))$,  and applying Proposition~\ref{prop:tildeM}, we conclude that $\|\Qb_i(\bmu(t))\| = O(\sigma_t)$. Thus, repeating the steps in~\eqref{eq:G_1_norm} - ~\eqref{eq:firstEst}, where we count for the additional term $\Qb_i(\bmu(t))$  we obtain:	
	\begin{align}
		&\E\{\|\bmu^i(t+1)-\zb^i(t)\|^2|\EuScript F_t\} \cr
		&\le\|\bmu^i(t)-\zb^i(t)\|^2 \cr
		&\quad- 2\gamma_t({\Mb}_i(\bmu(t))- {\Mb}_i(\zb(t)), \bmu^i(t)-\zb^i(t))\cr
		&\quad + 2\gamma_t\|\Qb(t)\|\|\bmu^i(t)-\zb^i(t)\|\cr
		& \quad+ \gamma^2(t)[\|{\Mb}_i(\bmu(t)) - {\Mb}_i(\zb(t))\|^2+\|\Qb(t)\|^2+\E\{\|\Rb_{i,j}(t)\|^2|\EuScript F_t\}].
	\end{align}
	Thus, using continuity of ${\Mb}$ and the fact that both $\bmu^i(t)$ and $\zb^i(t)$ belong to the compact set $\Ab$, as well as taking into account the result  of  Lemma~\ref{lem:Rsq}, we conclude that
	\begin{align}\label{eq:secondEst1}
		&\E\{\|\bmu^i(t+1)-\zb^i(t)\|^2|\EuScript F_t\} \cr
		&\le\|\bmu^i(t)-\zb^i(t)\|^2 \cr
		&\quad- 2\gamma_t({\Mb}_i(\bmu(t))- {\Mb}_i(\zb(t)), \bmu^i(t)-\zb^i(t))+ h_1(t),
	\end{align}
	where
	$h_1(t) = O(\g_t\sigma_t+\g_t^2).
	$
	Thus, due to Assumption~\ref{assum:CG_grad}, we conclude from~\eqref{eq:secondEst1} by summing up the inequalities over $i=1,\ldots,N$,
	\[\E\{\|\bmu(t+1)-\zb(t)\|^2|\EuScript F_t\}\le(1-\nu\gamma_t)\|\bmu(t)-\zb(t)\|^2+ Nh_1(t).\]
	Next, by using the inequality
	\[\|\bmu(t)-\zb(t)\|^2\le (1+\epsilon_t)\|\bmu(t)-\zb(t-1)\|^2 + (1+1/\epsilon_t)\|\zb(t-1) -\zb(t)\|^2,\]
	which holds for any $\epsilon_t>0$,
	we get
	\begin{align*}
		\E\{\|\bmu(t+1)-\zb(t)\|^2|\EuScript F_t\}&\le(1-\nu\gamma_t)(1+\epsilon_t)\|\bmu(t)-\zb(t-1)\|^2 \cr
		&\,  + (1-\nu\gamma_t)(1+1/\epsilon_t)\|\zb(t-1) -\zb(t)\|^2+h_1(t)\cr
		&\le\left(1-\frac{\nu\gamma_t}{2}\right)\|\bmu(t)-\zb(t-1)\|^2 \cr
		&\qquad+ O\left(\frac{|\rho_t-\rho_{t-1}|^2}{\gamma_t}+h_1(t)\right),
	\end{align*}
	where in the last inequality we chose $\epsilon_t = \frac{\nu\gamma_t}{2}$ and used Proposition~\ref{prop:diffz}.
	Thus, given the settings for the parameters, the definition of $h_1(t)$,  we conclude that
	\begin{align}\label{eq:LV2z}
		\E\{\|\bmu(t+1)-\zb(t)\|^2|\EuScript F_t\}\le\left(1-\frac{2}{t}\right)\|\bmu(t)-\zb(t-1)\|^2 + O\left(\frac{1}{t^{2}}\right).
	\end{align}
	Thus, \eqref{eq:LV2z} implies that $\bmu(t)$ converges to $\zb(t-1)$ almost surely (see Chung's Lemma in Appendix~\ref{app:Ch}). Moreover, due to Proposition~\ref{prop:z-mu}, the almost sure convergence of $\bmu(t)$ to $\ab^*$ takes place. Taking into account that $\bxi(t)\sim\EuScript N(\bmu(t),\sigma_t)$ and $\sigma_t\to 0$ as $t\to\infty$, we conclude that $\bxi(t)$ converges weakly to a Nash equilibrium $\ab^*$. Moreover, according to the Portmanteau Lemma \cite{portlem}, this convergence is also in probability.
	Next, by taking the full expectation of the both sides in~\eqref{eq:LV2}, we obtain
	\begin{align*}
		&\E[\|\bmu(t+1)-\zb(t)\|^2]\le\left(1-\frac{2}{t}\right)\E\|\bmu(t) - \zb(t-1)\|^2+ O\left(\frac{1}{t^{2}}\right).
	\end{align*}
	By applying the Chung's lemma (see Appendix~\ref{app:Ch}) to the inequality above, we conclude that
	\[\E[\|\bmu(t+1)-\zb(t)\|^2] = O\left(\frac{1}{t}\right).\]
		Finally,
	\begin{align}\label{eq:final}
		\E[\|\bmu(t+1)-\ab^*\|^2]  & \le  2\E[\|\bmu(t+1)-\zb(t)\|^2] +  2\|\zb(t)-\ab^*\|^2.
	\end{align}
	Thus, applying Proposition~\ref{prop:z-mu} with the given choice of the parameters, we obtain the result.
\end{proof}

\subsection{Discussion}
With two-point function evaluation, the proposed payoff-based procedure achieves the best convergence rate $O(1/T)$ in the class of stochastic gradient algorithms applied to optimization of smooth strongly convex functions, given bounded estimations of gradients at each time~\cite{Rakhlin2012MakingGD}. This improved rate is due to the fact that, in this case, expectation of the squared stochastic term $\Rb_{i,2}(t)$ in the proposed method~\eqref{eq:pbavmu} is upper bounded by some constant (see Lemma~\ref{lem:Rsq}). This rate is however not achievable in the case of one-point estimations of gradients since in that case we have unbounded expectation as over time ($\sigma_t\to 0$ as time runs), which requires a specific setting for the variance parameter $\s_t$ to be able to upper bound the corresponding stochastic term during the iterates. Such settings naturally lead to a slower convergence rate.

With one-point function evaluation, the proposed payoff-based procedure achieves the best convergence rate $O(1/\sqrt{T})$ in the class of stochastic gradient algorithms applied to optimization of smooth strongly convex functions, given zero-order function evaluation at each time~\cite{shamir2013complexity}. On the other hand, the work~\cite{bravo2018bandit} addresses convergence to Nash equilibria in strictly monotone games given payoff information with one-point gradient estimations. In particular, the authors mentioned the rate of their proposed algorithm in the case of the strongly monotone pseudo-gradients. The obtained rate was $O\left(\frac{1}{t^{1/3}}\right)$ which is worse than $O\left(\frac{1}{t^{1/2-\epsilon}}\right)$ presented in this paper. The reason is two different techniques used in the procedures' analysis. The work~\cite{bravo2018bandit} estimates the distances from the iterates to the actual solution and not to the solutions in games with mixed strategies (in our notations such solutions are $\bmu^*_t$ and $\by(t)$ for the action set regularized by $\rho_t$).  The latter approach allows for getting a tighter upper bound for the distance between the iterates and the Nash equilibrium in the original game by obtaining the rate of convergence to the mixed-strategy solutions as well as estimating the distance between them and the actual one (see Propositions~\ref{prop:distSol}, \ref{prop:y-mu}), and finally applying the triangle inequality (see~\eqref{eq:final} in the proof of Theorem~\ref{th:main}).

\section{Conclusion}
\label{sec:conclusion}
We showed that bandit learning of the Nash equilibrium in convex strongly monotone smooth games can be achieved with a rate of $O\left(\frac{1}{\sqrt{T}}\right)$ in the one-point feedback and $O\left(\frac{1}{T}\right)$ in the two-point feedback. These rates are the lowest known to our knowledge for the considered class of problems.

In future, we need to quantify the rate dependence on the problem dimension $d$ and number of players $N$. Furthermore, we may verify that  the mirror-descent class of algorithms achieve the same optimal rates. Finally, characterizing the optimal rates for a larger class of games, namely, those with non-strongly monotone pseudo-gradients remains an open problem.

\bibliographystyle{plain}
\bibliography{srtrMonGames_ref}

%%%%%%%%%%%%%%%%%%%%%%%%%%%%%%%%%%%%%%%%%%%%%%%%%%%%%%%%%%%%%%%%%%%%%%%%%%%%%%%%%%%%%%%%%
%%%%%%%%%%%%%%%%%%%%%%%%%%
%%%%%%%%%%%%%%%%APPENDIX%%%%%%%%%%%%%%%%%
%%%%%%%%%%%%%%%%%%%%%%%%%%
%%%%%%%%%%%%%%%%%%%%%%%%%%%%%%%%%%%%%%%%%%%%%%%%%%%%%%%%%%%%%%%%%%%%%%%%%%%%%%%%%%%%%%%%%

\appendix
\section{Auxiliary Results}
The following auxiliary lemma will be used in the proofs of some propositions and lemmas below.
\begin{lem}\label{lem:aux}
	Let some continuous function $f(\bx): \R^{Nd}\to\R$ be  such that $f(\bx)\ge 0$ for any $\bx\in\R^{Nd}$ and  $f(\bx) = O(\exp\{\|\bx\|^{\alpha}\})$ as $\|\bx\|\to\infty$, where $\alpha<2$. Let $\Ab$ be some compact subset of $\R^{Nd}$. Finally, let $p(\bx; \bmu, \sigma_t)$, $\bmu\in\Ab$, be the density function of the Gaussian vector $\bxi$ as defined in~\eqref{eq:densityfull} and $\tilde \bxi = \bmu+\theta(\bxi-\bmu)$ for some $\theta\in[0,1]$.
	Then there exists a constant $C>0$  such that $\E\{f(\tilde \bxi)\}\le 2^{Nd}C$.
	Moreover, if $\sigma_t\to 0$ as $t\to\infty$, then there exists $c>0$, which is independent on $d$ and $N$, such that $\E\{f(\tilde \bxi)\}\le c$ for all sufficiently large $t$.
\end{lem}
\begin{proof}
	Let us define $\tilde \bx = \bmu+\theta(\bx-\bmu)$ for any $\bx\in\R^{Nd}$. Then,
	\begin{align*}
	&\E\{f(\tilde \bxi)|\bxi\sim\EuScript N(\bmu,\sigma_t)\} = \int_{\R^{Nd}}f(\tilde \bx)p( \bx; \bmu, \sigma_t)d\bx\\
	&=\int_{\mathbb R^{Nd}}f(\tilde\bx)\prod_{i=1}^{N} \frac{1}{(\sqrt{2\pi}\sigma_t)^{d}}\exp\left\{-\sum_{k=1}^{d}\frac{(x^i_k-\mu^i_k)^2}{2\sigma_t^2}\right\}d\bx\cr
	& = 2^{Nd}\int_{\mathbb R^{Nd}}f(\tilde\bx)\prod_{i=1}^{N} \exp\left\{-\sum_{k=1}^{d}\frac{(x^i_k-\mu^i_k)^2}{4\sigma_t^2}\right\}\cr
	&\qquad\qquad \times \frac{1}{(\sqrt{2\pi}2\sigma_t)^{Nd}}\prod_{i=1}^{N} \exp\left\{-\sum_{k=1}^{d}\frac{(x^i_k-\mu^i_k)^2}{4\sigma_t^2}\right\}d\bx.
\end{align*}
According to the condition $f(\bx) = O(\exp\{\|\bx\|^{\alpha}\})$ as $\|\bx\|\to\infty$, where $\alpha<2$, there exists a constant  $C>0$  such that
\[f(\tilde\bx)\prod_{i=1}^{N} \exp\left\{-\sum_{k=1}^{d}\frac{(x^i_k-\mu^i_k)^2}{4\sigma_t^2}\right\}\le C\]
for any $\bx\in\R^{Nd}$ and $\bmu\in\Ab$.
Thus, taking into account that
\[\int_{\mathbb R^{Nd}}\frac{1}{(\sqrt{2\pi}2\sigma_t)^{Nd}}\prod_{i=1}^{N} \exp\left\{-\sum_{k=1}^{d}\frac{(x^i_k-\mu^i_k)^2}{4\sigma_t^2}\right\}d\bx = 1,\]
we conclude that $\E{f(\bxi)}\le 2^{Nd}C$.

Furthermore,  let us notice that
\begin{align}\label{eq:boundR}
	\E\{f(\tilde\bxi)\}& = \int_{\mathbb R^{Nd}}f(\tilde\bx) p(\bx;\bmu,\s_t)d\bx\cr
	& =\int_{\mathbb R^{Nd}}f(\tilde\bx)\prod_{k=1}^{N} \frac{1}{(\sqrt{2\pi}\sigma_{t})^{d}}\exp\left\{-\sum_{k=1}^{d}\frac{(x^i_k-\mu^i_k)^2}{2\sigma^2_{t}}\right\}d\bx\cr
	&\le c_0 + \int_{\mathbb R^{Nd}\setminus \boldsymbol{B}}f(\tilde\bx)\prod_{k=1}^{N} \frac{1}{(\sqrt{2\pi}\sigma_{t})^{d}}\exp\left\{-\sum_{k=1}^{d}\frac{(x^i_k-\mu^i_k)^2}{2\sigma^2_{t}}\right\}d\bx,
\end{align}
where $\boldsymbol{B}$ is some compact such that $\Ab\subset \boldsymbol{B}$ and we used boundedness of $f$ over $\boldsymbol{B}$ to conclude that there exists some constant $c_0$ such that
\[\int_{\boldsymbol{B}}f(\tilde\bx)\prod_{k=1}^{N} \frac{1}{(\sqrt{2\pi}\sigma_{t})^{d}}\exp\left\{-\sum_{k=1}^{d}\frac{(x^i_k-\mu^i_k)^2}{2\sigma^2_{t}}\right\}d\bx\le c_0.\]
   Next, we notice that, as $\sigma_t\to 0$, and due to the fact that $f(\bx) = O(\exp\{\|\bx\|^{\alpha}\})$ as $\|\bx\|\to\infty$, where $\alpha<2$,  for given $d$ and $N$, there exists some finite $T$ such that the following inequality holds for $t\ge T$ and any $\bx\in \mathbb R^{Nd}\setminus \boldsymbol{B}$:
	$$f(\tilde\bx)\prod_{k=1}^{N} \exp\left\{-\sum_{k=1}^{d}\frac{(x^i_k-\mu^i_k)^2}{4\sigma^2_{t}}\right\}\le \frac{1}{2^{Nd}}.$$
Hence, for $t\ge T$
\begin{align}\label{eq:boundRR}
	&\int_{\mathbb R^{Nd}\setminus \boldsymbol{B}}f(\tilde\bx)\prod_{k=1}^{N} \frac{1}{(\sqrt{2\pi}\sigma_{t})^{d}}\exp\left\{-\sum_{k=1}^{d}\frac{(x^i_k-\mu^i_k)^2}{2\sigma^2_{t}}\right\}d\bx \cr
	& =\int_{\mathbb R^{Nd}\setminus \boldsymbol{B}}f(\tilde\bx)\prod_{k=1}^{N} \exp\left\{-\sum_{k=1}^{d}\frac{(x^i_k-\mu^i_k)^2}{4\sigma^2_{t}}\right\} \frac{1}{(\sqrt{2\pi}\sigma_{t})^{Nd}}\prod_{k=1}^{N} \exp\left\{-\sum_{k=1}^{d}\frac{(x^i_k-\mu^i_k)^2}{4\sigma^2_{t}}\right\}d\bx\cr
	&\le \frac{1}{2^{Nd}}\int_{\mathbb R^{Nd}\setminus \boldsymbol{B}} \frac{1}{(\sqrt{2\pi}\sigma_{t})^{Nd}}\prod_{k=1}^{N} \exp\left\{-\sum_{k=1}^{d}\frac{(x^i_k-\mu^i_k)^2}{4\sigma^2_{t}}\right\}d\bx\cr
	&\le \int_{\mathbb R^{Nd}}\frac{1}{(\sqrt{2\pi}2\sigma_{t})^{Nd}}\prod_{k=1}^{N} \exp\left\{-\sum_{k=1}^{d}\frac{(x^i_k-\mu^i_k)^2}{4\sigma^2_{t}}\right\}d\bx = 1.
\end{align}
Thus, we conclude from~\eqref{eq:boundR} and~\eqref{eq:boundRR} that  $\E\{{J_i}^4(\bxi(t))\}\le c_0+1 =c$ in the case $\sigma_t\to 0$.
\end{proof}

We will use the H\"older's inequality:
\begin{align}\label{eq:HI}
	\E|XY|\le (\E(X^2))^{1/2}(\E(Y^2))^{1/2},
\end{align}
which holds for any random variables $X$ and $Y$ such that $\E(X^2)<\infty$ and $\E(Y^2)<\infty$.

\section{Proof of Lemma~\ref{lem:sample_grad}}\label{app:sample}
\begin{proof}
	First, we verify that the differentiation under the integral sign in
	\begin{align*}
		\tilde M{i,k} (\bmu,\sigma)=\frac{\partial {\int_{\mathbb R^{Nd}}J_i(\bx)p( \bx; \bmu, \sigma)d\bx}}{\partial \mu^i_k}, \mbox{ for $k\in[d]$},
	\end{align*}
	is justified at any $\bmu\in\Ab$. Then, it can then readily be verified that  the first equality in~\eqref{eq:deriv} holds, by taking the differentiation inside the integral.
	
	A sufficient condition for differentiation under the integral is that the integral of the formally differentiated function with respect to $\mu^i_k$  converges uniformly over the domain set of the parameter $\mu^i_k$, whereas the differentiated function is continuous (see \cite[Chapter 17]{zorich}). Continuity of the functions $J_i(\bx)p( \bx; \bmu, \sigma)$ and
	$\frac{\partial {J_i(\bx)p( \bx; \bmu, \sigma)}}{\partial \mu^i_k} = J_i(\bx)\frac{x^i_k -\mu^i_k}{\sigma^2}p( \bx; \bmu, \sigma)$
	follows from continuity of $J_i(\bx)$ and $p( \bx; \bmu, \sigma)$. We demonstrate uniform convergence of the integral
	$\int_{\R^{Nd}} J_i(\bx)\frac{x^i_k -\mu^i_k}{\sigma^2}p( \bx; \bmu, \sigma)d\bx$
	over the set $\Ab$. By making the substitution $\yb = \xb - \bmu$, we obtain
	\begin{align*}
		&\int_{\R^{Nd}} J_i(\bx)\frac{x^i_k -\mu^i_k}{\sigma^2}p( \bx; \bmu, \sigma)d\bx  \\
		&=\int_{\R^{Nd}} J_i(\by+\bmu)\frac{y^i_k}{\sigma^2}p( \by; 0, \sigma)d\by.
	\end{align*}
	Thus, taking into account Assumption~\ref{assum:infty} and the fact that $\by\in\Ab$, we conclude that there exists a finite constant $D>0$ such that for any $\yb\in\R^{Nd}$
	\begin{align*}
		|J_i(\by+\bmu)\frac{y^i_k}{\sigma^2}p( \by; 0, \sigma)|\le D \frac{|y^i_k|}{\sigma^2} \prod_{i=1}^{N} \frac{1}{(\sqrt{2\pi}\sigma)^{d}}\exp\left\{-\sum_{k=1}^{d}\frac{(y^i_k)^2}{4\sigma^2} \right\}.
	\end{align*}
	By noting that
	\begin{align*}
	 &\int_{\R^{Nd}} D \frac{|y^i_k|}{\sigma^2} \prod_{i=1}^{N} \frac{1}{(\sqrt{2\pi}\sigma)^{d}}\exp\left\{-\sum_{k=1}^{d}\frac{(y^i_k)^2}{4\sigma^2} \right\}d\yb = O(\sigma),
	\end{align*}
	we apply the Weierstrass criterion for uniform convergence of integrals to conclude the desired result.

	Next,
\begin{align}\label{eq:mixedstr_grad1}
	&\int_{\R^{Nd}}J_i(\bx)\frac{x_i^k - \mu^k_i}{\sigma^2}p(\bx;\bmu,\sigma)d\bx =-\frac{1}{(\sqrt{2\pi}\sigma)^{Nd}}\cr
	&\times\int_{\R^{(Nd-1)}}\left[\int_{x_i^k=-\infty}^{x_i^k=+\infty}J_i(\bx)d\left(\exp\left\{-\sum_{i=1}^{N}\sum_{k=1}^{d}\frac{(x_i^k-\mu^k_i)^2}
	{2\sigma^2}\right\}\right)\right]\cr
	&\qquad\times \exp\left\{-\sum_{j\ne i}^{N}\sum_{l\ne k}^d\frac{(x_j^l-\mu_j^l)^2}{2\sigma^2}\right\} d x^{-i,k}=\int_{\R^{Nd}}\frac{\partial J_i(\bx)}{\partial x_i^k}p(\bx;\bmu,\sigma)d\bx,
\end{align}
	where in the above, we use integration by parts and Assumption~\ref{assum:infty}.
\end{proof}

\section{Proof of Lemma~\ref{lem:Rsq}}\label{app:Rsq}
\begin{proof}
	From the definition of $\Rb_{i,j}$ in~\eqref{eq:Rterm} and its property given by \eqref{eq:mathexp2}, we notice that for $j=1$,
	\begin{align}
		\label{eq:Rineq1}
		\E\{\|\Rb_{i,1}(t)\|^2|\EuScript F_t\}&=\E\{\|\Rb_{i,1}(t)\|^2|\bxi(t)\sim \EuScript(\bmu(t),\sigma_t)\}\cr
		&\le\sum_{k=1}^{d}\int_{\mathbb R^{Nd}}{J_i}^2(\bx)\frac{(x^i_k - \mu^i_k(t))^2}{\sigma^4_t} p(\bx;\bmu(t),\s_t)d\bx
	\end{align}
	Next, we take into account Assumption~\ref{assum:infty} and apply Lemma~\ref{lem:aux} to conclude that there exists some finite constant $c_1>0$ such that
	\begin{align}\label{eq:1}
		\E\{{J_i}^4(\bxi(t))\}<c_1.
		\end{align}
	  Thus, we can use the H\"older's inequality~\eqref{eq:HI} to obtain
	\begin{align}\label{eq:boundR0}
		\int_{\mathbb R^{Nd}}&{J_i}^2(\bx)\frac{(x^i_k - \mu^i_k(t))^2}{\sigma^4_t} p(\bx;\bmu(t),\s_t)d\bx \cr
		&= \E\{{J_i}^2(\bxi(t))\frac{(\xi^i_k - \mu^i_k(t))^2}{\sigma^4_t}\}\le (\E\{{J_i}^4(\bxi(t))\})^{1/2}\left(\E\{\frac{(\xi^i_k - \mu^i_k(t))^4}{\sigma^8_t}\}\right)^{1/2}\cr
		& \le \frac{\sqrt{c_1}}{\sigma_t^2},
	\end{align}
	where in the last inequality we used \eqref{eq:1}. Hence,
	\[\E\{\|\Rb_{i,1}(t)\|^2|\EuScript F_t\}=O\left(\frac{d}{\sigma_t^2}\right).\]

	Now we consider $\Rb_{i,2}(t)$. According to the definition of $\Rb_{i,2}(t)$ in \eqref{eq:Rterm},
		\begin{align}
		\label{eq:Rineq2}
		\E\{\|\Rb_{i,2}(t)\|^2|\EuScript F_t\}&=\E\{\|\Rb_{i,2}(t)\|^2|\bxi(t)\sim \EuScript(\bmu(t),\sigma_t)\}\cr
		&\le\sum_{k=1}^{d}\int_{\mathbb R^{Nd}}({J_i}(\bx) - {J_i}(\bmu(t)))^2\frac{(x^i_k - \mu^i_k(t))^2}{\sigma^4_t} p(\bx;\bmu(t),\s_t)d\bx.
	\end{align}
  Next, due to the Taylor's expansion of $J_i$ around the point  $\bmu(t)$, we conclude that for any $\bx\in\R^{Nd}$
 \[{J_i}(\bx)={J_i}(\bmu(t)) +( \nabla J_i(\tilde{\bx}(t)),\bmu(t) - \bx)\]
where $\tilde{\bx}(t) = \bmu(t)+\theta(\bx-\bmu(t))$ for some $\theta\in[0,1]$. Thus, we obtain from~\eqref{eq:Rineq2}
 \begin{align}\label{eq:boundR2}
	\E\{\|\Rb_{i,2}(t)\|^2|\EuScript F_t\}&\le \int_{\mathbb R^{Nd}} \|\nabla J_i(\tilde{\bx}(t))\|^2\|\bmu(t) - \bx\|^2\left(\sum_{k=1}^{d}\frac{(x^i_k - \mu^i_k(t))^2}{\sigma^4_t}\right) p(\bx;\bmu(t),\s_t)d\bx\cr
	& =\frac{1}{\sigma^4_t} \E\left\{ \|\nabla J_i(\tilde{\bxi}(t))\|^2\|\bmu(t) - {\bxi}(t)\|^2\left(\sum_{k=1}^{d}{(\xi^i_k(t) - \mu^i_k(t))^2}\right)\right\},
\end{align}
where $\tilde{\bxi}(t) = \bmu(t)+\theta(\bxi(t)-\bmu(t))$. Next, we apply Lemma~\ref{lem:aux} (see Assumption~\ref{assum:infty}) to conclude existence of some constant $c_2>0$ such that
\begin{align}\label{eq:2}
	\E\|\nabla J_i(\tilde{\bxi}(t))\|^4\le c_2.
	\end{align}
Thus, applying the H\"older's inequality~\eqref{eq:HI} to~\eqref{eq:boundR2}, we get:
\begin{align*}
	\E\{\|\Rb_{i,2}(t)\|^2|\EuScript F_t\}\le& \sqrt{c_2}Nd^2.
\end{align*}
\end{proof}

\section{Proof of Lemma~\ref{lem:projTerm}}\label{app:proj}
\begin{proof}
	As $\ab(t) = \Proj_{\Ab}\bxi(t)$, we conclude that almost surely
	%Furthermore, because $\|\bxi(t)-\ab(t)\|\le \|\bxi(t)-\bxi\|$ for $\ab(t) = \Proj_{\Ab}\bxi(t)$. Thus, we conclude that
	\begin{align}\label{eq:projterm0}
		&\E\left\{|J_i(\ab(t))-J_i(\bxi(t))|\frac{\|\bxi^i(t) -\bmu^i(t)\|}{\sigma^2_t}\;|\EuScript F_t\right\} \cr
		&= \E\left\{|J_i(\ab(t))-J_i(\bxi(t))|\frac{\|\bxi^i(t) -\bmu^i(t)\|}{\sigma^2_t}| \bxi(t)\sim \EuScript N(\bmu(t),\s_t) \right\}\cr
		&\le (\E\left\{|J_i(\ab(t))-J_i(\bxi(t))|^2\right\})^{\frac{1}{2}}\left(\E\left\{\frac{\|\bxi(t)-\bmu(t)\|^2}{\sigma^4_t}\right\}\right)^{\frac{1}{2}}\cr
		& =\frac{\sqrt{Nd}}{\sigma_t}(\E\left\{|J_i(\ab(t))-J_i(\bxi(t))|^2\right\})^{\frac{1}{2}},
		%&\le \Pr\{\bxi(t)\in\R^{Nd}\setminus \Ab\} \E\left\{\frac{|J_i(\ab(t))-J_i(\bxi(t))| \|\bxi^i(t) -\bmu^i(t)\|}{\sigma^2_t}\right\}\cr
		%&\le \Pr\{\bxi(t)\in\R^{Nd}\setminus \Ab\} (\E\left\{|J_i(\ab(t))-J_i(\bxi(t))|^2\right\})^{\frac{1}{2}}\left(\E\left\{\frac{\|\bxi(t)-\bmu(t)\|^2}{\sigma^4_t}\right\}\right)^{\frac{1}{2}}\cr
		%&=\frac{\sqrt{Nd}}{\sigma_t}\Pr\{\bxi(t)\in\R^{Nd}\setminus \Ab\}(\E\left\{|J_i(\ab(t))-J_i(\bxi(t))|^2\right\})^{\frac{1}{2}},
		% k_0, \,\mbox{for some $k_0>0$},o
	\end{align}
	where, as before, we used the H\"older's inequality ($\E|XY|\le (\E(X^2))^{1/2}(\E(Y^2))^{1/2})$ and the notation $\E\{\cdot\} = \E\{\cdot |\bxi(t)\sim \EuScript(\bmu(t),\sigma_t)\}$.
	
	Next, we estimate the term $\E\left\{|J_i(\ab(t))-J_i(\bxi(t))|^2\right\}$ as follows.
	\begin{align}\label{eq:projterm}
		&\E\left\{|J_i(\ab(t))-J_i(\bxi(t))|^2\right\} \cr
		&=\int_{\Ab}|J_i(\Proj_{\Ab}\bx)-J_i(\bx)|^2p(\bx;\bmu(t),\sigma_t)d\bx\cr
		&\qquad\qquad + \int_{\R^{Nd}\setminus \Ab}|J_i(\Proj_{\Ab}\bx)-J_i(\bx)|^2p(\bx;\bmu(t),\sigma_t)d\bx\cr
		&=\int_{\R^{Nd}\setminus \Ab}|J_i(\Proj_{\Ab}\bx)-J_i(\bx)|^2p(\bx;\bmu(t),\sigma_t)d\bx\cr
		&\le 2\int_{\R^{Nd}\setminus \Ab}(J_i(\Proj_{\Ab}\bx))^2p(\bx;\bmu(t),\sigma_t)d\bx + 2\int_{\R^{Nd}\setminus \Ab}(J_i(\bx))^2p(\bx;\bmu(t),\sigma_t)d\bx\cr
		&\le 2K\Pr\{\bxi(t)\in\R^{Nd}\setminus \Ab\}) + 2\int_{\R^{Nd}\setminus \Ab}(J_i(\bx))^2p(\bx;\bmu(t),\sigma_t)d\bx,
	\end{align}
	where $O$-notation above is in terms of $t\to\infty$.
	The second equality above is due to the fact that $\Proj_{\Ab}\bx=\bx$ for any $\bx\in \Ab$. The first inequality was obtained by taking into account that $|J_i(\Proj_{\Ab}\bx)-J_i(\bx)|^2\le 2 (J_i(\Proj_{\Ab}\bx))^2 + 2(J_i(\bx))^2$, whereas the last one is due to the inequality $(J_i(\Proj_{\Ab}\bx))^2\le K$ for some constant $K$ and any $\bx$ (since $\Ab$ is compact and $J_i$ is continuous).

	% Analogously to \eqref{eq:boundR}, we can demonstrate that $\E\left\{|J_i(\ab(t))-J_i(\bxi(t))|^2\right\}$ is bounded.
	Thus, let us estimate $\Pr\{\bxi(t)\in\R^{Nd}\setminus \Ab\}$ above. The idea is that since $\bxi(t)$ is sampled from a Gaussian distribution with mean $\bmu(t)$, $\bxi(t)$ concentrates around its mean $\bmu(t)$ with high probability. Since the mean is projected onto a shrunk version of the set $\Ab$, namely, $(1-\rho_t) \Ab$, by appropriately tuning $\rho_t$ and the variance of the distribution $\sigma_t$ we can ensure that $\bxi(t)$ stays within the original feasible set with high probability.
	
	Let $\EuScript O_{\rho_t}(\bx) = \{\boldsymbol y\in\R^{Nd}| \|\boldsymbol y-\bx\|^2<{\rho}_t^2\}$ denote the $\rho_t$-neighborhood of the point $\bx\in\Ab$. Hence, $\sup_{\boldsymbol y\notin \EuScript O_{\rho_t}(\bx)}-\|\boldsymbol y-\bx\|^2=-{\rho}_t^2$ . Then, taking into account the fact that $\EuScript O_{\rho_t}(\bx)$ is contained in $\Ab$ and $\rho_t<1$, we obtain that for any $t$ and any bounded $\sigma>\sigma_t$:
	\begin{align}\label{eq:probterm}
		\Pr&\{\bxi(t)\in\R^{Nd}\setminus \Ab\}\le\Pr\{\bxi(t)\in\R^{Nd}\setminus \EuScript O_{\rho_t}(\bmu(t))\}\cr
		& = \int_{\boldsymbol y\notin \EuScript O_{\rho_t}(\bmu(t))}\frac{1}{(2\pi)^{Nd/2}\sigma_t^{Nd}}
		\exp\left\{-\frac{\|\boldsymbol y-\bmu(t)\|^2}{2\sigma_t^2}\right\}d\boldsymbol y \cr
		&=\int_{\boldsymbol y\notin \EuScript O_{\rho_t}(\bmu(t))}\exp\left\{-\|\boldsymbol y-\bmu(t)\|^2\left(\frac{1}{2\sigma_t^2} - \frac{1}{2\sigma^2}\right)\right\}\cr
		&\qquad\qquad\times\frac{\sigma^{Nd}}{\sigma_t^{Nd}}\frac{1}{(2\pi)^{{Nd}/2}\sigma^{Nd}}\exp\left\{-\frac{\|\boldsymbol y-\bmu(t)\|^2}{2\sigma^2}\right\}d\boldsymbol y\cr
		&\le \exp\left\{-\rho^2_t\left(\frac{1}{2\sigma_t^2} - \frac{1}{2\sigma^2}\right)\right\}\frac{\sigma^{Nd}}{\sigma_t^{Nd}}\cr
		&\quad\times\int_{\boldsymbol y\notin \EuScript O_{\rho_t}(\bmu(t))}\frac{1}{(2\pi)^{Nd/2}\sigma^{Nd}}\exp\left\{-  \frac{\|\boldsymbol y-\bmu(t)\|^2}{2\sigma^2}\right\}d\boldsymbol y\cr
		&\le k_1 \frac{e^{-\frac{{\rho}_t^2}{2\sigma_t^2}}}{\sigma_t^{Nd}},
	\end{align}
	for some finite $k_1>0$. The last inequality holds because
	\[\int_{\boldsymbol y\notin \EuScript O_{\rho_t}(\bmu(t))}\frac{1}{(2\pi)^{{Nd}/2}\sigma^{Nd}}\exp\left\{-  \frac{\|\boldsymbol y-\bmu(t)\|^2}{2\sigma^2}\right\}d\boldsymbol y\le 1\]
	and, thus, due to the diminishing $\rho_t$ there exists $0<k_1<\infty$:
	\[\int_{\boldsymbol y\notin \EuScript O_{\rho_t}(\bmu(t))}\frac{e^{ \frac{{\rho}_t^2}{2\sigma^2}}\sigma^{Nd}}{(2\pi)^{N/2}\sigma^{Nd}}\exp\left\{-  \frac{\|\boldsymbol y-\bmu(t)\|^2}{2\sigma^2}\right\}d\boldsymbol y\le k_1.\]
	
	Taking into account Assumption~\ref{assum:infty}, we conclude existence of some constant $C_1$ such that \[(J_i(\xb))^2\exp\left\{-\frac{\|\boldsymbol x-\bmu(t)\|^2 }{4\sigma_t^2}\right\}\le C_1\]
	for any $\xb$. Thus,
	\begin{align*}
		&\int_{\R^{Nd}\setminus \Ab}(J_i(\bx))^2p(\bx;\bmu(t),\sigma_t)d\bx\cr
		&\le C_1\int_{\boldsymbol y\in \R^{Nd}\setminus \Ab}\frac{1}{(2\pi)^{Nd/2}\sigma_t^{Nd}}
		\exp\left\{-\frac{\|\boldsymbol x-\bmu(t)\|^2 }{4\sigma_t^2}\right\}d\boldsymbol x\cr
		& \le C_1\int_{\boldsymbol x\notin \EuScript O_{\rho_t}(\bmu(t))}\frac{1}{(2\pi)^{Nd/2}\sigma_t^{Nd}}
		\exp\left\{-\frac{\|\boldsymbol x-\bmu(t)\|^2 }{4\sigma_t^2}\right\}d\boldsymbol x.
	\end{align*}
	Hence, analogously to \eqref{eq:probterm}, we have
	\begin{align}\label{eq:probterm1}
		&\int_{\R^{Nd}\setminus \Ab}(J_i(\bx))^2p(\bx;\bmu(t),\sigma_t)d\bx \cr
		& \le C_1\int_{\boldsymbol x\notin \EuScript O_{\rho_t}(\bmu(t))}\frac{1}{(2\pi)^{Nd/2}\sigma_t^{Nd}}
		\exp\left\{-\frac{\|\boldsymbol x-\bmu(t)\|^2 }{4\sigma_t^2}\right\}d\boldsymbol x\cr
		&=\int_{\boldsymbol x\notin \EuScript O_{\rho_t}(\bmu(t))}\exp\left\{-\|\boldsymbol x-\bmu(t)\|^2\left(\frac{1}{4\sigma_t^2} - \frac{1}{4\sigma^2}\right)\right\}\cr
		&\qquad\qquad\times\frac{(\sqrt{2}\sigma)^{Nd}}{\sigma_t^{Nd}}\frac{1}{(2\pi)^{{Nd}/2}(\sqrt{2}\sigma)^{Nd}}\exp\left\{-\frac{\|\boldsymbol x-\bmu(t)\|^2}{4\sigma^2}\right\}d\boldsymbol x\cr
		&\le \exp\left\{-\rho^2_t\left(\frac{1}{4\sigma_t^2} - \frac{1}{4\sigma^2}\right)\right\}\frac{(\sqrt{2}\sigma)^{Nd}}{\sigma_t^{Nd}}\cr
		&\quad\times\int_{\boldsymbol x\notin \EuScript O_{\rho_t}(\bmu(t))}\frac{1}{(2\pi)^{Nd/2}(\sqrt{2}\sigma)^{Nd}}\exp\left\{-  \frac{\|\boldsymbol x-\bmu(t)\|^2}{4\sigma^2}\right\}d\boldsymbol x\cr
		&\le k_2 \frac{e^{-\frac{{\rho}_t^2}{2\sigma_t^2}}}{\sigma_t^{Nd}}= O\left(\frac{e^{-\frac{{\rho}_t^2}{2\sigma_t^2}}}{\sigma_t^{Nd}}\right),
	\end{align}
	where the last inequality is due to the fact that
	\begin{align*}
		&\int_{\boldsymbol x\notin \EuScript O_{\rho_t}(\bmu(t))}\frac{1}{(2\pi)^{Nd/2}(\sqrt{2}\sigma)^{Nd}}\exp\left\{-  \frac{\|\boldsymbol x-\bmu(t)\|^2}{4\sigma^2}\right\}d\boldsymbol x\\
		&\le \int_{\R^{Nd}}\frac{1}{(2\pi)^{Nd/2}(\sqrt{2}\sigma)^{Nd}}\exp\left\{-  \frac{\|\boldsymbol x-\bmu(t)\|^2}{4\sigma^2}\right\}d\boldsymbol x=1.
	\end{align*}
\end{proof}

\section{Proof of Proposition~\ref{prop:tildeM}}\label{app:convex}
\begin{proof}
	Due to differentiability of $\Mb$ and its Taylor's expansion, for any $i\in[N]$, $k\in[d]$,
	\begin{align}\label{eq:Taylor}
			M_{i,k}(\bx)=M_{i,k}(\bmu) &+ (\nabla M_{i,k}(\tilde{\bx}), \bx-\bmu),
	\end{align}
where $\tilde{\bx} = \bmu+\theta(\bx-\bmu)$ for some $\theta\in[0,1]$.
Thus,
\begin{align}\label{eq:tildeM1}
&\left|\tilde M_{i,k}^{(t)} (\bmu)-\tilde M_{i,k}^{(t-1)}(\bmu)\right|  = \left|\tilde M_{i,k}^{(t)} (\bmu)-M_{i,k}(\bmu)+M_{i,k}(\bmu)-\tilde M_{i,k}^{(t-1)}(\bmu)\right|\cr
&=\left|\int_{\R^{Nd}}[M_{i,k}(\bx)-M_{i,k}(\bmu)]p(\bx;\bmu,\sigma_t)d\bx+\int_{\R^{Nd}}[M_{i,k}(\bmu)-M_{i,k}(\bx)]p(\bx;\bmu,\sigma_{t-1})d\bx\right|\cr
& = \big|\int_{\R^{Nd}}(\nabla M_{i,k}(\tilde{\bx}), \bx-\bmu)p(\bx;\bmu,\sigma_t)d\bx\cr
&\qquad-\int_{\R^{Nd}}(\nabla M_{i,k}(\tilde{\bx}), \bx-\bmu)p(\bx;\bmu,\sigma_{t-1})d\bx\big|,
\end{align}
where in the last equality we used~\eqref{eq:Taylor}. Moreover let us consider the function $f_i^k(\bmu,\sigma_t) =\tilde M_{i,k}^{(t)}(\bmu)  - M_{i,k}(\bmu)$ separately. Thus, according to ~\eqref{eq:Taylor},
	\begin{align*}
		f_i^k(\bmu,\sigma_t) =&\int_{\R^{Nd}}[M_{i,k}(\bx)-M_{i,k}(\bmu)]p(\bx;\bmu,\sigma_t)d\bx\cr
		=&\int_{\R^{Nd}}(\nabla M_{i,k}(\tilde{\bx}), \bx-\bmu)p(\bx;\bmu,\sigma_t)d\bx\cr
		=&\E\{(\nabla M_{i,k}(\tilde{\bxi}), \bxi-\bmu)\},
	\end{align*}
   where $\bxi$ is the Gaussian vector with the density function $p(\bx;\bmu,\sigma_t)$ and $\tilde{\bxi} = \bmu+\theta(\bxi-\bmu)$.
	Thus, using Cauchy-Schwarz inequality, applicable due to Lemma~\ref{lem:aux} (see Assumption~\ref{assum:infty}), we obtain
	\begin{align}\label{eq:f}
		|f_i^k(\bmu,\sigma_t)| \le&\E\{\|\nabla M_{i,k}(\tilde{\bxi})\|\| \bxi-\bmu\|\}\cr
		\le&\left(\E\{\|\nabla M_{i,k}(\tilde{\bxi})\|^2\}\right)^{1/2}\left(\E\{\| \bxi-\bmu\|^2\}\right)^{1/2} = O(\sqrt{Nd}\sigma_t),
	\end{align}
 where the second inequality is due to~\eqref{eq:HI} and the last equality follows from Lemma~\ref{lem:aux} and the fact that $\E \|\bxi-\bmu\|^2 = \tat{Nd\sigma_t^2}$.
	Hence,
	\[\|\tilde{\Mb}^{(t)} - \Mb(\bmu)\| = O(\tat{(Nd)^{\frac{3}{2}}\sigma_t}).\]
	Moreover, taking into account \eqref{eq:tildeM1}, \eqref{eq:f}, and given $\boldsymbol f =(f_1^1,\ldots,f_i^k,\ldots f_N^d)$, we conclude that
	\[\|\tilde{\Mb}^{(t)} (\bmu)-\tilde{\Mb}^{(t-1)}(\bmu)\| =\|\boldsymbol f(\bmu,\sigma_t)-\boldsymbol f(\bmu,\sigma_{t-1})\| = O((Nd)^{\frac{3}{2}|\sigma_t - \sigma_{t-1}|}).\]
\end{proof}

\section{Proof of Proposition~\ref{prop:strMon}}\label{app:strMon}
\begin{proof}
	We start by characterizing the elements of the matrix $\frac{\partial \tilde {\Mb}^{(t)}(\bmu)}{\partial \bmu}$ with $\bmu\in\Ab$.
Let us consider the element $\frac{\partial \tilde M_{i,k}^{(t)}(\bmu)}{\partial \mu_j^l}$, $k,l\in[d]$, $j\in[N]$. We have
\[\frac{\partial \tilde M_{i,k}^{(t)}(\bmu)}{\partial \mu_j^l} = \frac{\partial \int_{\R^{Nd}}M_{i,k}(\bx)p(\bx;\bmu,\sigma_t)d\bx}{\partial \mu_j^l}.\]
Analogously to the proof of Lemma~\ref{lem:sample_grad},  we can demonstrate that in the expression above the differentiation under the integral sign is justified for all $k,l\in[d]$, $j\in[N]$. Hence, differentiability of $\tilde{\Mb}^{(t)}$ and Lipschitz continuity over compact sets follow.
Thus,
\begin{align}\label{eq:mainnu}
	\frac{\partial \tilde M_{i,k}^{(t)}(\bmu)}{\partial \mu_j^l} &=  \int_{\R^{Nd}}\frac{\partial M_{i,k}(\bx)p(\bx;\bmu,\sigma_t)d\bx}{\partial \mu_j^l} = -\int_{\R^{Nd}}(x_i^l - \mu_i^l) M_{i,k}(\bx)p(\bx;\bmu,\sigma_t)d\bx\cr
	& = -\frac{1}{(\sqrt{2\pi}\sigma_t)^{Nd}}\int_{\R^{Nd-1}}\left[\int_{x_i^l=-\infty}^{x_i^l=+\infty}M_{i,k}(\bx)
	d\left(\exp\left\{-\frac{(x_i^l-\mu_i^l)^2}{2\sigma_t^2}\right\}\right)\right]\cr
	&\qquad\qquad\times\exp\left\{\frac{-\|\bx_{-i}-\bmu_{-i}\|^2}{2\sigma_t^2}\right\} d \bx_{-i}\cr
	&=\int_{\R^{Nd}}\frac{\partial M_{i,k}(\bx)}{\partial x_i^l}p(\bx;\bmu,\sigma_t)d\bx \cr
	&= \int_{\R^{Nd}}\left(\frac{\partial M_{i,k}(\bx)}{\partial x_i^l} -\frac{\partial M_{i,k}(\bmu)}{\partial \mu_i^l}\right) p(\bx;\bmu,\sigma_t)d\bx + \frac{\partial M_{i,k}(\bmu)}{\partial \mu_i^l},
\end{align}
where in the third equality above we used integration by parts.
Let us denote the function $\frac{\partial M_{i,k}(\bx)}{\partial x_i^l}$ by $g_{i,k,l}(\bx)$, i.e. $g_{i,k,l} =\frac{\partial M_{i,k}(\bx)}{\partial x_i^l}$.
 If Assumption~\ref{assum:Lipschitz}.\ref{num:11} holds, then we can apply the Taylor's expansion to the function $g_{i,k,l}(\cdot)$ at the point $\bmu$ to conclude that
 \begin{align*}
 g_{i,k,l}(\bx) - g_{i,k,l}(\bmu) = (\nabla g_{i,k,l}(\tilde \bx), \bx-\bmu),
 \end{align*}
where $\tilde \bx = \bmu+\theta(\bx-\bmu)$ for some $\theta\in[0,1]$.
Hence,
 \begin{align}\label{eq:HIto1}
	\int_{\R^{Nd}}\left(g_{i,k,l}(\bx) - g_{i,k,l}(\bmu)\right) p(\bx;\bmu,\sigma_t)d\bx&\ge- \int_{\R^{Nd}}\|\nabla g_{i,k,l}(\tilde \bx)\|\|\bx-\bmu\| p(\bx;\bmu,\sigma_t)d\bx\cr
	&=-\E\{\|\nabla g_{i,k,l}(\tilde \bxi)\|\|\bxi-\bmu\|\},
\end{align}
where $\bxi$ is the Gaussian vector with the density function $p(\bx;\bmu,\sigma_t)$ and $\tilde \bxi = \bmu+\theta(\bxi-\bmu)$.
Next, according to Lemma~\ref{lem:aux}, there exists $K_1$ such that
\[\E\{\|\nabla g_{i,k,l}(\tilde \bxi)\|^2\}\le K_1^2.\]
Moreover,
\begin{align*}
	\E\{\|\bxi-\bmu\|^2\}= \int_{\R^{Nd}}\|\bx-\bmu\| p(\bx;\bmu,\sigma_t)d\bx & = \sum_{i=1}^{N}\sum_{k=1}^{d}\int_{\R^{Nd}}\|x^k_i-\mu^k_i\|^2 p(\bx;\bmu_t,\sigma_t)d\bx  = Nd\sigma^2_t.
\end{align*}
Thus, using~\eqref{eq:HI} in~\eqref{eq:HIto1}, we obtain
 \begin{align}\label{eq:HIto2}
	\int_{\R^{Nd}}\left(g_{i,k,l}(\bx) - g_{i,k,l}(\bmu)\right) p(\bx;\bmu,\sigma_t)d\bx&\ge- K_1\sqrt{Nd}\sigma_t.
\end{align}
	Further, by combining the elements $\frac{\partial \tilde M_{i,k}^{(t)}(\bmu)}{\partial \mu_j^l}$ above in the matrix $\frac{\partial \tilde {\Mb}^{(t)}(\bmu)}{\partial \bmu}$ and taking into account Assumption~\ref{assum:CG_grad} and, thus, the inequality $\frac{\partial  \Mb(\bmu)}{\partial \bmu}  \succcurlyeq \nu I_{Nd}$ with $I_{Nd}$ being the identity $Nd\times Nd$-matrix, we obtain, according to~\eqref{eq:mainnu}, that
	\begin{align*}
		\frac{\partial \tilde {\Mb}^{(t)}(\bmu)}{\partial \bmu} & \succcurlyeq -K_1Nd\sigma_t \mathbf{1}_{Nd} + \nu I_{Nd},
	\end{align*}
	where $\mathbf{1}_{Nd}$  in the $Nd\times Nd$-matrix with all elements equal to 1 and  of the same size.
	Thus, by taking $\sigma_t \le \frac{\nu}{2K_1N^2d^2}$, we obtain the strict diagonally dominant symmetric matrix $\frac{\partial \tilde {\Mb}^{(t)}(\bmu)}{\partial \bmu}$ such that
	\begin{align*}
		\frac{\partial \tilde {\Mb}^{(t)}(\bmu)}{\partial \bmu} & \succcurlyeq  \frac{\nu}{2} I_{Nd}.
	\end{align*}
On the other hand, if Assumption~\ref{assum:Lipschitz}.\ref{num:12} holds, we obtain the last two inequalities above by directly applying the Lipschitz property of $\Mb$ to~\eqref{eq:mainnu} and replacing the constant $K_1$ by the Lipschitz constant $K$.
\end{proof}

	\section{Proof of Proposition~\ref{prop:distSol}}\label{app:distSol}
	\begin{proof}
		We use the well-known fact that $y^*\in SOL(Y,T)$ if and only if
		\[y^* = \Proj_{Y}[y^*-\beta T(y^*)]\]
		for  any  $\beta>0$.
		Thus, using the fact above and Theorem~\ref{th:VINE} we conclude that
		\begin{align*}
			\ab^*&=\Proj_{\Ab}[\ab^* - \beta \Mb(\ab^*)],\cr
			\bmu_t^*&=\Proj_{\Ab}[\bmu_t^* - \beta \tilde{\Mb}^{(t)}(\bmu_t^*),]
		\end{align*}
		where $\beta$ is any positive constant. Thus, using the non-expansion of the projection operator, we obtain
		\begin{align*}
			\|\bmu_t^*-\ab^*\|^2 & \le \|\bmu_t^* - \ab^* - \beta (\tilde{\Mb}^{(t)}(\bmu_t^*) - \Mb(\ab^*))\|^2 \cr
			& =\|\bmu_t^*-\ab^*\|^2 - 2 \beta (\tilde{\Mb}^{(t)}(\bmu_t^*) - \Mb(\ab^*),\bmu_t^*-\ab^*) + \beta^2 \|\Mb(\ab^*) - \tilde{\Mb}^{(t)}(\bmu_t^*)\|^2.
		\end{align*}
		By taking into account that, due to Assumptions~\ref{assum:CG_grad} - \ref{assum:Lipschitz} and \eqref{eq:MvsTildeM} in Proposition~\ref{prop:tildeM},
		\begin{align*}
			(\tilde{\Mb}^{(t)}(\bmu_t^*) - \Mb(\ab^*),\bmu_t^*-\ab^*) &= (\tilde{\Mb}^{(t)}(\bmu_t^*) -\Mb(\bmu_t^*)+\Mb(\bmu_t^*)-\Mb(\ab^*),\bmu_t^*-\ab^*)\cr
			&\ge - \|\tilde{\Mb}^{(t)}(\bmu_t^*)-\Mb(\bmu_t^*) \|\|\bmu_t^*-\ab^*\| + \nu\|\bmu_t^*-\ab^*\|^2 \cr
			& \ge  - k_1\sigma_t\|\bmu_t^*-\ab^*\| + \nu\|\bmu_t^*-\ab^*\|^2
		\end{align*}
		for some constant $k_1>0$. Moreover, due to compactness of $\Ab$ and continuity of $\Mb$, the mapping $\Mb$ is Lipschitz continuous over $\Ab$ with some constant $L$. Hence,
		\begin{align*}
			\|\Mb(\ab^*) - \tilde{\Mb}^{(t)}(\bmu_t^*)\|^2 & \le  2\|\Mb(\ab^*) - \Mb(\bmu_t^*)\|^2 + 2\|\Mb(\bmu_t^*)- \tilde{\Mb}^{(t)}(\bmu_t^*)\|^2\cr
			& \le 2L^2\|\bmu_t^*-\ab^*\|^2 + 2k_1^2\sigma_t^{2}.
		\end{align*}
		Thus, combining the inequalities above we obtain
		\begin{align*}
			\|\bmu_t^*-\ab^*\|^2 & \le \|\bmu_t^*-\ab^*\|^2 + 2\beta k_1\sigma_t\|\bmu_t^*-\ab^*\| - 2\nu\beta\|\bmu_t^*-\ab^*\|^2+2L^2\beta^2\|\bmu_t^*-\ab^*\|^2\cr
			&\qquad\qquad + 2k_1^2\beta^2\sigma_t^2.
		\end{align*}
		Hence,
		\[(\nu - \beta L^2)\|\bmu_t^*-\ab^*\|^2 -  k_1\sigma_t\|\bmu_t^*-\ab^*\|  - \beta k_1^2\sigma_t^{2} \le 0. \]
		By taking $\beta = \frac{\nu}{2L^2}$ we obtain the following inequality:
		\[\frac{\nu}{2}\|\bmu_t^*-\ab^*\|^2 -  k_1\sigma_t\|\bmu_t^*-\ab^*\|  - \frac{\nu k_1^2\sigma_t^{2}}{2L^2}  \le 0,\]
		which implies
		\[\|\bmu_t^*-\ab^*\|\le \frac{k_1\sigma_t(1+\sqrt{1+\frac{\nu^2}{L^2}})}{\nu}.\]
		 Thus, the result follows.
	\end{proof}
	
	\section{Proof of Proposition~\ref{prop:tVSt-1}}\label{app:tVSt-1}
	\begin{proof}
		We focus here on such $t$ for which Proposition~\ref{prop:strMon} holds, namely the mapping $\tilde{\Mb}^{(t)}$ is strongly monotone over $\R^{Nd}$ with the constant $\nu/2$.
		
		According to the definition of $\bmu^*_t$,
		\[\bmu^*_t = \Proj_{\Ab}[\bmu^*_t -\beta\tilde{\Mb}^{(t)}(\bmu^*_t) ],\]
		where $\beta$ is any positive constant.
		Thus, using the non-expansion of the projection operator, we obtain
		\begin{align*}
			\|\bmu^*_t-\bmu^*_{t-1}\| &\le \|\bmu^*_t-\bmu^*_{t-1} +\beta(\tilde{\Mb}^{(t-1)}(\bmu^*_{t-1}) -\tilde{\Mb}^{(t)}(\bmu^*_{t}))\|\cr
			&\le \|\bmu^*_t-\bmu^*_{t-1} +\beta(\tilde{\Mb}^{(t)}(\bmu^*_{t-1}) -\tilde{\Mb}^{(t)}(\bmu^*_{t}))\| \cr
			&\qquad\qquad+ \beta\|\tilde{\Mb}^{(t-1)}(\bmu^*_{t-1}) -\tilde{\Mb}^{(t)}(\bmu^*_{t-1})\|.
		\end{align*}
		Moreover,  the mapping $\tilde{\Mb}^{(t)}$ is Lipschitz continuous over $\Ab$ due to compactness of the latter (see Proposition~\ref{prop:strMon}). Note that, according again to Proposition~\ref{prop:strMon}, $L_{(t)}\ge\frac{\nu}{2}$. Thus,
		\begin{align}\label{eq:eq1}
			&\|\bmu^*_t-\bmu^*_{t-1} +\beta(\tilde{\Mb}^{(t)}(\bmu^*_{t-1}) -\tilde{\Mb}^{(t)}(\bmu^*_{t}))\|^2\cr
			&= \|\bmu^*_t-\bmu^*_{t-1}\|^2 +2\beta(\tilde{\Mb}^{(t)}(\bmu^*_{t-1}) -\tilde{\Mb}^{(t)}(\bmu^*_{t}),\bmu^*_t-\bmu^*_{t-1})\cr
			&\qquad\qquad\qquad+\beta^2\|\tilde{\Mb}^{(t)}(\bmu^*_{t-1}) -\tilde{\Mb}^{(t)}(\bmu^*_{t})\|^2\cr
			&\le(1 - \beta\nu)\|\bmu^*_t-\bmu^*_{t-1}\|^2 +\beta^2 L^2_{(t)}\|\bmu^*_t-\bmu^*_{t-1}\|^2 \cr
			& = (1 - \beta\nu +\beta^2 L^2_{(t)})\|\bmu^*_t-\bmu^*_{t-1}\|^2  = (1-\frac{\nu^2}{4L_{(t)}^2})\|\bmu^*_t-\bmu^*_{t-1}\|^2,
		\end{align}
		where we substituted  $\beta = \frac{\nu}{2L^2_{(t)}}$.
		Thus,
		\begin{align*}
			&\|\bmu^*_t-\bmu^*_{t-1}\|\le \sqrt{1-\frac{\nu^2}{4L^2_{(t)}}}\|\bmu^*_t-\bmu^*_{t-1}\| + \frac{\nu}{2L^2_{(t)}}\|\tilde{\Mb}^{(t-1)}(\bmu^*_{t-1}) -\tilde{\Mb}^{(t)}(\bmu^*_{t-1})\|,
		\end{align*}
		and, hence, by taking into account Proposition~\ref{prop:tildeM}~\eqref{eq:tildeM}, we obtain
		\begin{align*}
			\|\bmu^*_t-\bmu^*_{t-1}\|\le \sqrt{1-\frac{\nu^2}{4L_{(t)}^2}}\|\bmu^*_t-\bmu^*_{t-1}\| &+ \frac{\nu}{2L^2_{(t)}}O(|\sigma_t-\sigma_{t-1}|).
		\end{align*}
		which implies the result.
	\end{proof}
	
	\section{Proof of Proposition~\ref{prop:diff}}\label{app:diff}
	\begin{proof}
		We focus here on such $t$ for which Proposition~\ref{prop:strMon} holds, namely the mapping $\tilde{\Mb}^{(t)}$ is strongly monotone over $\R^{Nd}$ with the constant $\nu/2$.
		
		According to the definition of $\yb(t)$,
		\[\yb(t) = \Proj_{(1-\rho_t)\Ab}[\yb(t) -\alpha\tilde{\Mb}^{(t)}(\yb(t)) ],\]
		where $\alpha$ is any positive constant.
		Thus, using the triangle inequality, Lemma~\ref{lem:A_r} (see Appendix~\ref{app:aux}), and the non-expansion property of the projection operator, we obtain
		\begin{align*}
			&\|\yb(t)-\yb(t-1)\| \cr
			&\le\|\Proj_{(1-\rho_t)\Ab}[\yb(t) -\tilde{\Mb}^{(t)}(\yb(t)) ] - \Proj_{(1-\rho_t)\Ab}[\yb(t-1) -\tilde{\Mb}^{(t-1)}(\yb(t-1)) ]\| \cr
			&+\|\Proj_{(1-\rho_t)\Ab}[\yb(t-1) -\tilde{\Mb}^{(t-1)}(\yb(t-1))] - \Proj_{(1-\rho_{t-1})\Ab}[\yb(t-1) -\tilde{\Mb}^{(t-1)}(\yb(t-1)) ]\|\|\cr
			&\le\|\Proj_{(1-\rho_t)\Ab}[\yb(t) -\tilde{\Mb}^{(t)}(\yb(t)) ] - \Proj_{(1-\rho_t)\Ab}[\yb(t-1) -\tilde{\Mb}^{(t-1)}(\yb(t-1)) ]\| \cr
			&\qquad\qquad\qquad+ O(|\rho_t-\rho_{t-1}|)\cr
			& \le \|\yb(t)-\yb(t-1) +\alpha(\tilde{\Mb}^{(t-1)}(\yb(t-1)) -\tilde{\Mb}^{(t)}(\yb(t)))\|+ O(|\rho_t-\rho_{t-1}|).
		\end{align*}
		We proceed with estimating $\|\yb(t)-\yb(t-1) +\alpha(\tilde{\Mb}^{(t-1)}(\yb(t-1)) -\tilde{\Mb}^{(t)}(\yb(t)))\|$.
		\begin{align*}
			&\|\yb(t)-\yb(t-1) +\alpha(\tilde{\Mb}^{(t-1)}(\yb(t-1)) -\tilde{\Mb}^{(t)}(\yb(t)))\|\cr
			&\le \|\yb(t)-\yb(t-1) +\alpha(\tilde{\Mb}^{(t)}(\yb(t-1)) -\tilde{\Mb}^{(t)}(\yb(t)))\| \cr
			&\qquad\qquad+ \alpha\|\tilde{\Mb}^{(t-1)}(\yb(t-1)) -\tilde{\Mb}^{(t)}(\yb(t-1))\|.
		\end{align*}
		Moreover,  the mapping $\tilde{\Mb}^{(t)}$ is Lipschitz continuous over $\Ab$ with some constant $L_{(t)}$ (see Proposition~\ref{prop:strMon}). Note that, according again to Proposition~\ref{prop:strMon}, $L_{(t)}\ge\frac{\nu}{2}$. Thus,
		\begin{align}\label{eq:eq1}
			&\|\yb(t)-\yb(t-1) +\alpha(\tilde{\Mb}^{(t)}(\yb(t-1)) -\tilde{\Mb}^{(t)}(\yb(t)))\|^2\cr
			&= \|\yb(t)-\yb(t-1)\|^2 -2\alpha(\tilde{\Mb}^{(t)}(\yb(t-1)) -\tilde{\Mb}^{(t)}(\yb(t)),\yb(t-1)-\yb(t))\cr
			&\qquad\qquad\qquad+\alpha^2\|\tilde{\Mb}^{(t)}(\yb(t-1)) -\tilde{\Mb}^{(t)}(\yb(t))\|^2\cr
			&\le(1 - \alpha\nu)\|\yb(t)-\yb(t-1)\|^2 +\alpha^2 L^2_{(t)}\|\yb(t)-\yb(t-1)\|^2 \cr
			& = (1 - \alpha\nu +\alpha^2 L^2_{(t)})\|\yb(t)-\yb(t-1)\|^2  = (1-\frac{\nu^2}{4L_{(t)}^2})\|\yb(t)-\yb(t-1)\|^2,
		\end{align}
		where substituted  $\alpha = \frac{\nu}{2L^2_{(t)}}$.
		Thus,
		\begin{align*}
			&\|\yb(t)-\yb(t-1) +\alpha(\tilde{\Mb}^{(t-1)}(\yb(t-1)) -\tilde{\Mb}^{(t)}(\yb(t)))\|\cr
			&\le \sqrt{1-\frac{\nu^2}{4L^2_{(t)}}}\|\yb(t)-\yb(t-1)\| + \frac{\nu}{2L^2_{(t)}}\|\tilde{\Mb}^{(t-1)}(\yb(t-1)) -\tilde{\Mb}^{(t)}(\yb(t-1))\|,
		\end{align*}
		and, hence, by taking into account \eqref{eq:tildeM} in Proposition~\ref{prop:tildeM}, we obtain
		\begin{align*}
			\|\yb(t)-\yb(t-1)\|\le \sqrt{1-\frac{\nu^2}{4L_{(t)}^2}}\|\yb(t)-\yb(t-1)\| &+ \frac{\nu}{2L^2_{(t)}}O(|\sigma_t-\sigma_{t-1}|)\cr
			&+ O(|\rho_t-\rho_{t-1}|).
		\end{align*}
		which implies the result.
	\end{proof}
	
	\section{Proof of Proposition~\ref{prop:y-mu}}\label{app:y-mu}
	\begin{proof}
		We focus here on such $t$ for which Proposition~\ref{prop:strMon} holds, namely the mapping $\tilde{\Mb}^{(t)}$ is strongly monotone over $\R^{Nd}$ with the constant $\nu/2$.
		
		According to the definition of $\bmu_t^*$ and $\yb(t)$, we have
		\begin{align*}
			\yb(t) & = \Proj_{(1-\rho_t)\Ab}[\yb(t) -\alpha\tilde{\Mb}^{(t)}(\yb(t))],\cr
			\bmu^*_t & =\Proj_{\Ab}[\bmu^*_t -\alpha\tilde{\Mb}^{(t)}(\bmu^*_t) ],
		\end{align*}
		given any $\alpha>0$. Thus, taking into account Lemma~\ref{lem:A_r} (see Appendix~\ref{app:aux}), we obtain
		\begin{align}\label{eq:eq2}
			\|\yb(t)-\bmu^*_t\|&\le \|\Proj_{(1-\rho_t)\Ab}[\yb(t) -\alpha\tilde{\Mb}^{(t)}(\yb(t)) ] - \Proj_{(1-\rho_t)\Ab}[\bmu^*_t -\alpha\tilde{\Mb}^{(t)}(\bmu^*_t) ]\|\cr
			& \qquad+ \|\Proj_{(1-\rho_t)\Ab}[\bmu^*_t -\alpha\tilde{\Mb}^{(t)}(\bmu^*_t) ] - \Proj_{\Ab}[\bmu^*_t -\alpha\tilde{\Mb}^{(t)}(\bmu^*_t) ]\|\|\cr
			&\le \|\yb(t)-\bmu^*_t -\alpha(\tilde{\Mb}^{(t)}(\yb(t)) -\alpha\tilde{\Mb}^{(t)}(\bmu^*_t))\| + O(\rho_t).
		\end{align}
		We proceed with estimating $\|\yb(t)-\bmu^*_t -\alpha(\tilde{\Mb}^{(t)}(\yb(t)) -\alpha\tilde{\Mb}^{(t)}(\bmu^*_t))\|$.
		Analogously to~\eqref{eq:eq1}, we get for $\alpha = \frac{\nu}{2L^2_{(t)}}$:
		\begin{align*}
			&\|\yb(t)-\bmu^*_t -\alpha(\tilde{\Mb}^{(t)}(\yb(t)) -\alpha\tilde{\Mb}^{(t)}(\bmu^*_t))\|^2\cr
			& \le (1-\frac{\nu^2}{4L_{(t)}^2})\|\yb(t)-\bmu^*_t\|^2,
		\end{align*}
		where $L_{(t)}$ is, as before, the Lipschitz constant of the mapping $\tilde{\Mb}^{(t)}$. Applying the result above to~\eqref{eq:eq2}, we obtain the result.
	\end{proof}

	\section{Properties of the shrinked set}
	\label{app:aux}
	\begin{lem}\label{lem:A_r}
		For any $\xb\in\R^{Nd}$ the following holds:
		\[\|\Proj_{(1-\rho_{t-1})\Ab}\xb - \Proj_{(1-r_{t})\Ab}\xb\| = O(|\rho_{t-1}-\rho_t|),\]
		\[\|\Proj_{\Ab}\xb - \Proj_{(1-r_{t})\Ab}\xb\| = O(|\rho_t|).\]
	\end{lem}
	\begin{proof}
		%Note that by definition, for any $\xb\in\R^{Nd}$
		%\begin{align}\label{eq:A_r}
		%  \Proj_{(1-\rho_{t-1})\Ab}\xb = \Proj_{(1-\rho_{t-1})\Ab}\{\Proj_{(1-r_{t})\Ab}\xb\}.
		%\end{align}
		Without loss of generality,  assume  $\xb\notin (1-\rho_{t-1})\Ab$ (otherwise, $\|\Proj_{(1-\rho_{t-1})\Ab}\xb - \Proj_{(1-r_{t})\Ab}\xb\| = 0$).
		Due to convexity of the set $\Ab\subset \R^{Nd}$ there exists a convex function $g:\R^{Nd}\to\R$ such that $\Ab = \{\xb: g(\xb)\le 0\}$, so that $(1-\rho_t)\Ab = \{\xb: g(\xb)\le -\rho_t\}$ for any $t$. Moreover, define $\xb': = \Proj_{(1-r_{t})\Ab}\xb,$ and observe that $ \Proj_{(1-\rho_{t-1})\Ab}\xb = \Proj_{(1-\rho_{t-1})\Ab}\xb'$. Thus, we have $\|\Proj_{(1-\rho_{t-1})\Ab}\xb - \Proj_{(1-r_{t})\Ab}\xb\|= d,$ where $d$ is the optimal cost corresponding to the following optimization problem
		%where $d$ is defined below using $\xb' = \Proj_{(1-r_{t})\Ab}\xb$,
		\begin{align*}
			d := \min_{\yb} &\|\yb - \xb'\|, \\
			\mbox{s.t. } &g(\yb) = -\rho_{t-1}.
		\end{align*}
		The optimization problem has a solution $\yb^*$ for which the gradient of the corresponding Lagrangian is zero, namely
		\[\frac{(\yb^* - \xb')}{\|\yb^* - \xb'\|} + \lambda \nabla g(\yb^*) = \boldsymbol 0,\]
		where $\lambda>0$ is the dual multiplier of the problem under consideration. Due to Assumption~\ref{assum:convex} and the choice of $\rho_t$, that guarantees nonempty interior of $(1-\rho_t)\Ab$ for all $t$, the Slator's condition for the constraints $g(\xb)\le -\rho_t$ holds for all $t$. Hence, for any $\xb\in\R^{Nd}$ there exists a constant $\Lambda>0$ such that $\lambda<\Lambda$ (see \cite{ConvOpt}).
		Thus, we conclude that
		\[\nabla g(\yb^*) = -\frac{(\yb^* - \xb')}{ \lambda \|\yb^* - \xb'\|}.\]
		Next, due to convexity of the function $g$,
		\begin{align*}
			g(\xb')&\ge g(\yb^*) + (\nabla g(\yb^*), \xb'-\yb^*) \\
			& = -\rho_{t-1} + \frac{\|\yb^* - \xb'\|^2}{ \lambda \|\yb^* - \xb'\|}\ge -\rho_{t-1} + \frac{\|\yb^* - \xb'\|}{ \Lambda }.
		\end{align*}
		Thus, taking into account that $g(\xb') \le -r_{t}$, we obtain
		\[d =\|\yb^* - \xb'\|\le \Lambda (\rho_{t-1}-\rho_t) = O(|\rho_{t-1}-\rho_t |).\]
		Analogously one can demonstrate that for any $\xb\in\R^{Nd}$
		\[\|\Proj_{\Ab}\xb - \Proj_{(1-r_{t})\Ab}\xb\| = O(|\rho_t|).\]
	\end{proof}
	
	\section{The Chung's Lemma (Lemma 4 in Chapter 2.2. \cite{Polyak})}\label{app:Ch}
	\begin{lem}
		Let $u_k\ge 0$ and
		\[u_{k+1}\le\left(1-\frac{c}{k}\right)u_k + \frac{d}{k^{1+p}}, \quad d,p,c>0.\]
		Then
		\begin{align*}
			u_k\le d(c-p)^{-1}k^{-p} + o(k^{-p}), \quad &\mbox{if $c>p$},\cr
			u_k = O(k^{-c}\ln k), \quad &\mbox{if $c=p$}, \cr
			u_k = O(k^{-c}), \quad &\mbox{if $c<p$}.
		\end{align*}
	\end{lem}

	%
	% \section*{Acknowledgments}
	% We would like to acknowledge the assistance of volunteers in putting
	% together this example manuscript and supplement.

\end{document}